%% file: 2211.06442
\IfFileExists{filecontents.sty}{\RequirePackage{filecontents}}{}

\input{preamble}

\begin{document}
\input{abs}
\maketitle
\tableofcontents
\input{sec1-intro}

\input{sec2-prelims}

\input{sec3-construct}

\input{sec4-top-invar}
\input{sec5-examples}

\newpage
\input{sik.bbl}

\end{document}

%% file: preamble.tex
\documentclass[12pt]{amsart}
\usepackage[letterpaper,margin=1.25in]{geometry}
\usepackage[utf8]{inputenc}
\usepackage[T1]{fontenc}

\usepackage{
  amssymb,
  mathrsfs,
  enumitem,
  xcolor,
  tikz,
}

\usepackage{hyperref}
\hypersetup{
  pdftitle={Symplectic instanton knot homology},
  pdfauthor={David G. White},
  pdfborder={0 0 1},
  linkbordercolor=red,
  urlbordercolor=blue,
}

\usepackage{fancyhdr}
\allowdisplaybreaks

\makeatletter
\def\sectionmark#1{%
  \markboth {\MakeUppercase{%
      \ifnum \c@secnumdepth >\z@
	\thesection\quad
      \fi
#1}}{}}%

\makeatother

\DeclareMathOperator{\Hom}{Hom}
\DeclareMathOperator{\sik}{SIK}

\theoremstyle{plain}
\newtheorem{athm}{Theorem}
\renewcommand{\theathm}{\Alph{athm}}
\newtheorem{thm}{Theorem}[subsection]
\newtheorem*{thm*}{Theorem}
\newtheorem{cor}[thm]{Corollary}
\newtheorem*{cor*}{Corollary}
\newtheorem{lem}[thm]{Lemma}
\newtheorem*{lem*}{Lemma}
\newtheorem{prop}[thm]{Proposition}
\newtheorem*{prop*}{Proposition}
\newtheorem{conj}[thm]{Conjecture}
\newtheorem*{conj*}{Conjecture}

\theoremstyle{definition}
\newtheorem{defn}[thm]{Definition}
\newtheorem*{defn*}{Definition}
\newtheorem{exm}[thm]{Example}
\newtheorem*{exm*}{Example}

\newtheorem*{exr*}{Exercise}

\newtheorem*{prob*}{Problem}

\theoremstyle{remark}

\newtheorem*{obs*}{Observation}
\newtheorem{rem}[thm]{Remark}
\newtheorem*{rem*}{Remark}

\newtheorem*{note*}{Note}

\title{Symplectic instanton knot homology}
\author{David G. White}
\email{dgwhite2@ncsu.edu}
\date{\today}

%% file: abs.tex
\begin{abstract}
There have been a number of constructions of Lagrangian Floer homology 
invariants for $3$-manifolds defined in terms of symplectic character varieties 
arising from Heegaard splittings.  With the aim of establishing an Atiyah-Floer 
counterpart of Kronheimer and Mrowka's singular instanton homology, we 
generalize one of these, due to H. Horton, to produce a Lagrangian Floer 
invariant of a knot or link $K \subset Y$ in a closed, oriented $3$-manifold, 
which we call \emph{symplectic instanton knot homology} ($\mathrm{SIK}$).  We 
use a multi-pointed Heegaard diagram to parametrize the gluing together of a 
pair of handlebodies with properly embedded, trivial arcs to form $(Y, K)$.  
This specifies a pair of Lagrangian embeddings in the traceless 
$\mathrm{SU}(2)$-character variety of a multiply punctured Heegaard surface, 
and we show that this has a well-defined Lagrangian Floer homology.  Portions 
of the proof of its invariance are special cases of Wehrheim and Woodward's 
results on the quilted Floer homology associated to compositions of so-called 
elementary tangles, while others generalize their work to certain 
non-elementary tangles.
\end{abstract}

%% file: sec1-intro.tex
\section{Introduction}
We associate to a knot or link $K$ in a $3$-manifold $Y$ a Lagrangian Floer 
homology invariant $\mathrm{SIK}(Y, K)$, called \emph{symplectic instanton knot 
homology}.  By associating to the pair $(Y, K)$ a suitable multi-pointed 
Heegaard diagram, we obtain
\begin{itemize}
\item a smooth symplectic manifold in the form of the traceless 
$\mathrm{SU}(2)$-character variety of a closed, oriented punctured surface;
\item a pair of embedded monotone Lagrangian submanifolds parametrized by the 
Heegaard diagram.
\end{itemize}
By leveraging technical results of Oh \cite{Oh-stfh2} for the monotone case, we 
proceed to demonstrate that this gives rise to a well-defined Lagrangian Floer 
homology:

\begin{athm}[Floer homology] \label{thm:\theathm}
Let $Y$ be a connected, closed, oriented $3$-manifold and $K \subset Y$ a knot 
or link.  A choice of multi-pointed Heegaard diagram $\mathcal H$ compatible 
with the pair $(Y, K)$ determines a pair of Lagrangian submanifolds $L_\alpha, 
L_\beta$ in the traceless $\mathrm{SU}(2)$-character variety of the punctured 
Heegaard surface.  Their Lagrangian Floer homology
\begin{equation} \label{eq:sik-diagr-def}
\mathrm{SIK}(\mathcal H) := \mathrm{HF}(L_\alpha, L_\beta)
\end{equation}
is well defined with coefficients in $\mathbb Z$.
\end{athm}

The Heegaard diagram is constructed similarly as for knot and link Floer 
homology (see \cite{OS-holom_disk_knot_invar, Ras-floer_hom_knot_compl, 
OS-holom_disk_link_invar}), having $k$ pairs of marked points, considered as 
punctures, on the genus-$g$ Heegaard surface $\Sigma_g$ where it meets $K$ 
(i.e. $K$ is in $k$-bridge position).  However, a free basepoint is also taken 
and thickened into a triple of marked points.  Our symplectic manifold is the 
traceless $\mathrm{SU}(2)$-character variety of $\Sigma_g$ punctured at these 
marked points, denoted, $\mathscr R_{g,2k+3}$, with the extra triple of points 
serving to exclude reducible representations and thus ensure the character 
variety's smoothness.  The attaching curves of $\mathcal H$ determine two 
Lagrangian submanifolds $L_\alpha, L_\beta \subset \mathscr R_{g,2k+3}$, each 
of which is an embedding of the traceless character variety of one of the 
handlebodies in the Heegaard splitting, relative to $k$ properly embedded, 
trivial arcs where it meets $K$, along with three additional arcs forming a 
``tripod graph''.  We proceed to show that $(L_\alpha, L_\beta)$ is a monotone 
and relatively spin Lagrangian pair and prove Theorem \ref{thm:A} via results 
of Oh regarding such pairs.

It is well known that any two Heegaard diagrams of the kind utilized here are 
related by a finite sequence of certain operations, called Heegaard moves.  We 
study the impact of each of these moves on the Floer homology 
\eqref{eq:sik-diagr-def} using Wehrheim and Woodward's quilted Floer homology 
\cite{WW-qltd_floer_cohom} and find that any such move produces canonically 
isomorphic Floer homologies.  As a consequence we conclude that the isomorphism 
class of the Lagrangian Floer homologies \eqref{eq:sik-diagr-def} is an 
invariant of the pair $(Y, K)$ :
\begin{athm}[Topological invariance] \label{thm:\theathm}
Let $Y$ and $K$ be as in Theorem \ref{thm:A}, and suppose that $\mathcal H$ and 
$\mathcal H'$ are two multi-pointed Heegaard diagrams compatible with the pair 
$(Y, K)$.  Then the induced Lagrangian Floer homologies $\mathrm{SIK}(\mathcal 
H)$ and $\mathrm{SIK}(\mathcal H')$ are canonically isomorphic as relatively 
graded abelian groups.
\end{athm}
In light of Theorem \ref{thm:B}, we define the symplectic instanton knot 
homology as
\begin{equation} \label{eq:sik-def-intro}
\mathrm{SIK}(Y, K) := \text{ the isomorphism class of } \mathrm{SIK}(\mathcal H),
\end{equation}
for any choice of Heegaard diagram $\mathcal H$ compatible with $(Y, K)$.

\subsection{Motivation and antecedents} \label{sec:motiv-ante}

The principal goal of this construction is to carry out an Atiyah-Floer-type 
exercise, viz. to form a Lagrangian Floer counterpart of Kronheimer and 
Mrowka's \emph{singular instanton knot homology} \cite{KM-knot_hom_grp_inst, 
KM-kh_unknot_detector}.  The latter comes in two variants: the undreduced 
$I^\sharp (Y, K)$ and the reduced $I^\natural (Y, K)$.  By taking a certain 
auxiliary object $\theta$, called a theta graph, alongside $K$ as we construct 
a compatible Heegaard diagram, we obtain a Floer complex with generators which 
are identified via holonomy---at least before perturbations on either 
side---with those of $I^\sharp$.  The theta graph can be thought of as 
analogous to the Hopf link, with an arc connecting its components, employed in 
defining $I^\sharp$.  The computations in Section \ref{sec:examples} show that 
$\mathrm{SIK}$ agrees with $I^\sharp$ for the empty link and for the unknot in 
the $3$-sphere (see \cite[Lemma 
8.3]{KM-kh_unknot_detector}), and we hypothesize that the theories agree in 
  general, as predicted by the Atiyah-Floer conjecture:

\begin{conj}
Let $\mathcal H$ be a compatible diagram for $(Y, K)$.  Then there is an 
isomorphism $\mathrm{SIK}(\mathcal H) \cong I^\sharp(Y, K)$ as $\mathbb 
Z/2$-graded abelian groups.
\end{conj}

An Atiyah-Floer-type construction has already been undertaken for the reduced 
variant $I^\natural$ by Hedden, Herald and Kirk in the form of the their 
\emph{pillowcase homology} \cite{HHK-pillowcase1, HHK-pillowcase2}.  Their 
construction decomposes $(Y, K)$ along a Conway sphere, producing an 
intersection of immersed Lagrangian curves in the pillowcase.  A nice feature 
of this theory is the provision of explicit holonomy perturbations by which the 
authors obtain a minimal generating set for singular instanton homology.  Such 
concrete perturbation methods are currently absent from the theory we develop 
here, and whether similar holonomy perturbations can be applied in the context 
of $\mathrm{SIK}$ is an interesting question for study.

Our construction generalizes one of symplectic instanton homology for closed 
$3$-manifolds due to Horton \cite{Hor-sih_trless_xvar}, denoted 
$\mathrm{SI}(Y)$.  Indeed, $\mathrm{SIK}(Y, \emptyset) = \mathrm{SI}(Y)$, where 
$\emptyset$ denotes the empty link.  In the first instance, Horton's invariant 
is defined via Heegaard diagrams similarly as done here.  An appealing aspect 
of his theory, however, is that it admits an alternative definition as the 
quilted Floer homology \cite{WW-qltd_floer_cohom} induced by a Cerf 
decomposition of the $3$-manifold $Y$.  It is by means of this characterization 
that Horton is able to demonstrate both a K\"unneth principle for connected 
sums and a surgery exact triangle.  There is a natural analogue of this in the 
form of the quilted Floer homology associated to compositions of tangles.  A 
quilted Floer theory for tangles with trivial $3$-manifold topology has already 
been worked out in unpublished work of Wehrheim and Woodward 
\cite{WW-floer_fld_thry_tngl}.  We make a slight generalization of some of 
their results in proving the topological invariance of $\mathrm{SIK}$, with the 
expectation that a further generalization can be carried out, and that the 
formal properties of Horton's $\mathrm{SI}$ can thereby be extended to 
$\mathrm{SIK}$.  In particular, the example of (Heegaard) knot Floer homology 
$HFK$ \cite{OS-holom_disk_knot_invar, Ras-floer_hom_knot_compl} suggests that 
the surgery exact triangle for $\mathrm{SI}$ may be leveraged to demonstrate 
that $\mathrm{SIK}$ obeys a skein exact sequence.  Such a property is to be 
expected, given that instant knot homology, along with $HFK$, categorifies the 
classical Alexander polynomial.

\subsection{Organization of the paper} \label{sec:org}
In Section \ref{sec:prelim} we review material necessary for the proof of 
Theorem A, specifically the subjects of monotone Lagrangian Floer theory 
(Section \ref{sec:monot-lagr}) and of traceless character varieties of 
punctured surfaces (Section \ref{sec:xvar-surf}).  The construction of our 
invariant $\mathrm{SIK}$ is carried out in Section \ref{sec:constr}, wherein 
Theorem \ref{thm:A} is fully established.  Section \ref{sec:top-invar} 
addresses topological invariance, proving Theorem \ref{thm:B}.  In 
\ref{sec:examples} we compute $\mathrm{SIK}$ for the empty link and for unlinks 
in the $3$-sphere.

\subsection{Acknowledgements} \label{sec:ack}
The author wishes to express a debt of gratitude to his advisor, Tye Lidman, 
whose guidance and patience over the course of this and other projects has been 
invaluable.  The author was supported by NSF grants DMS-1709702 and 
DMS-2105469.  In addition part of this material is based upon work supported by 
National Science Foundation under Grant No. DMS-1928930, while the author was a 
Program Associate at the Simons Laufer Mathematical Sciences Institute 
(formerly \textsc{msri}) in Berkeley, California, during the Fall 2022 
semester.

%% file: sec2-prelims.tex
\section{Preliminaries} \label{sec:prelim}

\subsection{Monotone Lagrangian Floer theory} \label{sec:monot-lagr}

In \cite{Flo-morse_thry_lagr_inters} Floer introduced a homology theory 
associated to a pair of Lagrangian submanifolds $(L_0, L_1)$ in a symplectic 
manifold $(M, \omega)$, invariant up to Hamiltonian isotopy of the Lagrangians, 
called their \emph{Lagrangian Floer} (or \emph{Lagrangian intersection}) 
\emph{homology}.  In this intial formulation, strong constraints are placed on 
the topology of $M$, but subsequent elaborations of the theory have relaxed 
these considerably, at the expense of more or less complexity in the 
definition.  A broad treatment of the subject can be found in 
\cite{FO3-anom_obstr1, FO3-anom_obstr2}.  An especially convenient situation 
occurs when $(L_0, L_1)$ is a \emph{monotone} pair.  By work of Oh 
\cite{Oh-floer_cohom_lagr_int_i, Oh-floer_cohom_lagr_int_ii}, it is relatively 
straightforward to show under this condition that the obstructions to a 
well-defined Lagrangian Floer homology are averted.

The construction below produces a monotone Lagrangian pair, and our definition 
of symplectic instanton knot homology will rely on an easy corollary (Corollary 
\ref{cor:hf-gw-sympl}) of Oh's results.  We recall the key results below, 
preceded by a brief exposition to establish some notation and important notions 
such as relative spin structures appearing in the suite.  A thorough treatment 
can be found in Oh's two volumes \cite{Oh-stfh1, Oh-stfh2}, the second of which 
is the principal source for the material of this subsection.  Overviews of 
relevant subjects are given in e.g. \cite{CdS-sympl_geom} and 
\cite{Ped-qk_view_lagr_floer_hom}.

Let $(M, \omega)$ be a symplectic manifold.  For simplicity we will assume that 
$M$ is closed and connected.  Write $\mathcal J_\omega$ for the space of 
$\omega$-compatible almost complex structures on $M$, and let $c_1(M) \in 
H^2(M; \mathbb Z)$ be the first Chern class of $M$, defined as the first Chern 
class $c_1(TM, J)$ of the tangent bundle for any choice of $J \in \mathcal 
J_\omega$, being that the latter space is contractible.  Recall that $M$ is 
called monotone if $[\omega] \in H^2(M; \mathbb R)$ is proportional to 
$c_1(M)$, as made precise in the following:
\begin{defn}[Monotone symplectic manifold]
\label{def:monot-mfld}
A symplectic manifold $(M, \omega)$ is called ($\tau$-)\textbf{monotone} if 
there exists a positive real number $\tau$, called the \textbf{monotonicity 
constant}, such that
\begin{equation} \label{eq:monot-mfld}
[\omega] = \tau \cdot c_1(M).
\end{equation}
\end{defn}

\begin{rem} \label{rem:hurewicz}
Both $c_1(M), [\omega] : H_2(M) \to \mathbb Z, \mathbb R$ pull back to 
$\pi_2(M)$ via the Hurewicz homomorphism, and we use the same notation for 
these pullbacks.  Furthermore, because $\omega\vert_L \equiv 0$ by definition, 
there is a well-defined restriction of $[\omega]$ to $\pi_2(M, L)$.
\end{rem}

Now let $L \subset M$ be a Lagrangian submanifold. There is an associated map 
$\mu_L : \pi_2(M, L) \to \mathbb Z$ called the \emph{Maslov homomorphism}, and 
$L$ is called monotone when this and $[\omega]$ are proportional; more 
precisely:
\begin{defn}[Monotone Lagrangian] \label{def:monot-lagr}
A Lagrangian submanifold $L$ in a symplectic manifold $(M, \omega)$ is called 
\textbf{monotone} if there exists a positive real number $\kappa$ such that
\begin{equation} \label{eq:monot-lagr}
2[\omega]\vert_{\pi_2(M, L)} = \kappa \cdot \mu_L.
\end{equation}
\end{defn}

In general, if a Lagrangian $L \subset M$ is monotone, then $M$ must itself be 
monotone and their monotonicity constants equal, i.e. $\kappa = \tau$ (cf. 
\cite[Remark 2.3]{Oh-floer_cohom_lagr_int_i}).  There is a converse for simply 
connected Lagrangians:

\begin{prop}[] \label{prop:sc-monot}
If $M$ is a $\tau$-monotone symplectic manifold, then any compact simply 
connected Lagrangian $L \subset M$ is monotone with monotonicity constant 
$\tau$.
\end{prop}
\begin{proof}
One may cap off any $u : (D^2, \partial D^2) \to (M, L)$ with a disk in $L$ to 
form a sphere $v \subset M$.  By \cite[Appendix C]{MS-jholom_sympl_top} the 
Maslov index is given in this case in terms of the first Chern class by 
$\mu_L(u) = 2c_1(M)(v)$, which by the monotonicity of $M$ and the fact that 
$\omega\vert_L = 0$ is equal to $2\tau^{-1} [\omega](u)$, from which the result 
follows.
\end{proof}

Denote the positive generator of $\operatorname{im}(c_1(M)) \subset \mathbb Z$ 
by $N_M$, called the \emph{minimal Chern number}, and that of 
$\operatorname{im}(\mu_L) \subset \mathbb Z$ by $N_L$, called the \emph{minimal 
Maslov number}.  It is evident from the proof of Proposition 
\ref{prop:sc-monot} that
\begin{equation} \label{eq:maslov2chern}
N_L = 2 N_M
\end{equation}
for every compact simply connected Lagrangian $L \subset M$.

Along with monotonicity, a nice property of Lagrangians, and one which will be 
possessed by the ones constructed in the next section, is the following:

\begin{defn}[Relatively spin] \label{def:rel-spin}
A tuple $(L_0, \dotsc, L_m)$ of Lagrangian submanifolds in a symplectic 
manifold $M$ is said to be relatively spin if there exists a class $b \in 
H^2(M; \mathbb Z/2)$, called the \emph{background class}, which restricts to 
the second Stiefel-Whitney class of each Lagrangian: $b\vert_{L_i} = w_2(L_i)$ 
for $i = 0, \dotsc, m$.

Correspondingly, a \emph{relative spin structure} on $(L_0, \dotsc, L_m)$ is a 
choice of background class $b \in H^2(M; \mathbb Z/2)$ and of spin structures 
on $TL_i \oplus E$ where $E \to M$ is an orientable rank-$2$ vector bundle with 
$w_2(E) = b$.
\end{defn}

\begin{rem} \label{rem:spin}
In particular, if the Lagrangians of the above definition are all spin 
manifolds, then they are clearly relatively spin, as one may take $b = 0$ (cf.  
\cite[Sections 15.6.2-3]{Oh-stfh2}).
\end{rem}

When this additional structure is present, a choice of orientations of the 
Lagrangians canonically induces a coherent system of orientations upon the 
various moduli spaces of pseudo-holomorphic disks with boundaries on the 
Lagrangians.  As a result enumerative invariants derived from these moduli 
spaces may be counted with signs, taking values in $\mathbb Z$ rather than 
$\mathbb Z/2$; in particular, the Floer homology of a Lagrangian pair may be 
taken with integer coefficients.

We next establish one such enumerative invariant of monotone Lagrangians which 
will be important for our purposes, following \cite[Sec. 16.3]{Oh-stfh2}.  Let 
$L \subset M$ be a monotone Lagrangian with $N_L \ge 2$, and fix $\beta \in 
\pi_2(M, L)$ with $\mu_L(\beta) = 2$.  Given $J \in \mathcal J_\omega$, denote 
by $\mathcal M_1(L, \beta, J)$ the moduli space of $J$-holomorphic disks $u : 
(\mathbb D, \partial\mathbb D) \to (M, L)$ with $[u(\mathbb D)] = \beta$, 
modulo conformal automorphisms of the domain fixing $1 \in \mathbb D \subset 
\mathbb C$.  Consider the evaluation map
\begin{align}
ev_\beta : \mathcal M_1(L, \beta, J) & \to L \\
[u] & \mapsto u(1).
\end{align}
For any $x \in L$, one may choose $J$ generically such that $M_1(L, \beta, J)$ 
is a smooth compact manifold of dimension equal to that of $L$, and $ev_\beta$ 
is transverse to $x$.  Consequently there is a well-defined count
\begin{equation} \label{eq:disk-count}
\Phi(L, x) = \sum_{\beta \in \pi_2(M, L)} \# ev_\beta^{-1}(x),
\end{equation}
taken in $\mathbb Z/2$ in general, and in $\mathbb Z$ if $L$ is relatively 
spin.  Moreover, this count depends only on $[x] \in \pi_0(L)$; therefore, when 
$L$ is connected, we write simply $\Phi(L)$, called the \emph{disk number} of 
$L$.

\begin{lem}[] \label{lem:disk-num}
Let $L \subset M$ be a connected monotone Lagrangian and $\phi : M \to M$ a 
symplectomorphism.  Then $\Phi(\phi(L)) = \Phi(L)$.
\end{lem}
\begin{proof}
For a suitably generic choice of almost complex structure $J$, the 
symplectomorphism $\phi$ induces a diffeomorphism $\mathcal M_1(L, \beta, J) 
\approx \mathcal M_1(\phi(L), \phi_*\beta, \phi_*J)$.  Moreover, if $L$ is 
relatively spin and the disk numbers taken in $\mathbb Z$, one can choose 
relative spin structures compatibly, so that the induced diffeomorphism 
preserves the orientation systems on the moduli spaces.  Finally, $\phi_* : 
\pi_2(M, L) \to \pi_2(M, \phi(L))$ is an isomorphism, hence the sums 
\eqref{eq:disk-count} will be equal.
\end{proof}


The Lagrangians which will result from our construction below will turn out to 
have $N_L = 2$.  Oh shows in \cite[Sec. 16.4]{Oh-stfh2} that, for a pair of 
monotone Lagrangians $(L_0, L_1)$ with $N_L = 2$, the only potential 
obstruction to the Floer boundary operator's squaring to zero is so-called disk 
bubbling along the Lagrangians.  The following theorem provides conditions 
under which this is ruled out.

\begin{thm}[{\cite[Cor. 16.4.8]{Oh-stfh2}}]
\label{thm:hf-gw-invar}
Let $(L_0, L_1)$ be a pair of connected, compact monotone Lagrangian 
submanifolds in a symplectic manifold $(M, \omega)$ having minimal Maslov 
number $N_L = 2$.  For $i = 0, 1$, let $a_i$ be the positive generator of the 
group $\left\{ \omega(\beta) \mid \beta \in \pi_2(M, L_i) \right\}$, and let 
$\Phi(L_i)$ be the disk number of $L_i$.  Then $\partial \circ \partial = 0$ if 
and only if $\Phi(L_0) = \Phi(L_1)$ and $a_0 = a_1$.
\end{thm}

When the conditions of the forgoing theorem hold, the Lagrangian Floer homology 
$\mathrm{HF}(L_0, L_1)$ is well-defined.  As explained in \cite[Section 
15.6]{Oh-stfh2}, if $L_0, L_1$ are relatively spin, then the homology may be 
   taken with coefficients in $\mathbb Z$.

We will be able to apply Theorem \ref{thm:hf-gw-invar} in our case by virtue of 
the following special case:
\begin{cor} \label{cor:hf-gw-sympl}
Let $(L_0, L_1)$ be as in Theorem \ref{thm:hf-gw-invar}.  Suppose that there 
exists a symplectomorphism $\phi : (M, \omega) \to (M, \omega)$ such that 
$\phi(L_0) = L_1$.  Then the Lagrangian Floer homology $\mathrm{HF}(L_0, L_1)$ 
is well defined---with coefficients in $\mathbb Z/2$ generally and in $\mathbb 
Z$ if the pair $(L_0, L_1)$ is relatively spin.
\end{cor}
\begin{proof}
Since $\phi$ induces an isomorphism on $\pi_2(M, L_0) \to \pi_2(M, L_1)$ and 
preserves $\omega$, we have $a_0 = a_1$.  That $\Phi(L_0) = \Phi(L_1)$ follows 
from Lemma \ref{lem:disk-num}.
\end{proof}

\subsection{Traceless character varieties and punctured surfaces} 
\label{sec:xvar-surf}

Let $M$ be a manifold with or without boundary and $G$ be a compact Lie group.  
When $M$ is connected, its \emph{$G$-character variety} is the space $\chi_G(M) 
:= \Hom(\pi_1(M), G)/G$, where $G$ acts on representations by conjugation.  
When $M$ is disconnected with connected components $M_1, \dotsc, M_n$, we 
define its character variety as the cartesian product $\Pi_{i=1}^n 
\chi_G(M_i)$.  Now let $T \subset M$ be a properly embedded compact submanifold 
of codimension-$2$ having connected components $T_1, \dotsc, T_m$.  For each 
connected component $T_i$, one can choose a meridian $\mu_i$, the class of 
which is unique up to conjugacy in $\pi_1(M \setminus T)$.  We obtain a 
subspace of the relative character variety $\chi_G(M \setminus T)$ by 
specifying conjugacy classes $\mathcal C_1, \dotsc, \mathcal C_m$ in $G$ and 
requiring that representations take $[\mu_i]$ into $\mathcal C_i$, which we 
write as
\begin{equation} \label{eq:rel-xvar}
\chi_{G; \mathcal C_1, \dotsc, \mathcal C_m}(M, T) := \left\{ [\rho] \in 
\chi_G(M \setminus T) \mid \rho([\mu_i]) \in \mathcal C_i \text{ for } i = 1, 
\dotsc, m \right\}.
\end{equation}
Here we will concern ourselves only with the case where all of the $\mathcal 
C_i$ are equal to some fixed conjugacy class $\mathcal C \subset G$, for which 
we write just $\chi_{G, \mathcal C}(M, T)$.

The types of pairs $(M, T)$ with which we will work are \emph{punctured 
surfaces} and \emph{tangles}.  By a punctured surface we specifically mean that 
$M$ is a closed oriented surface, and $T$ is a finite collection of marked 
points, or punctures, on $M$ oriented by a map $\varepsilon : T \to \left\{ \pm 
1 \right\}$.  By a tangle we mean that $M$ is a connected compact oriented 
  $3$-manifold (possibly with boundary) and $T$ a properly embedded compact 
oriented $1$-manifold.  For our purposes the crucial result regarding the
character varieties of such spaces is the following:
\begin{thm}[{\cite[Cor. 4.5]{CHK-tngl_xvar}}]
\label{thm:lagr-immers}
Let $G$ be a compact Lie group and $\mathcal C \subset G$ an 
inversion-invariant conjugacy class.  Let $(Y, T)$ be a tangle, as defined 
above, where $Y$ has the possibly disconnected boundary surface $\partial Y = 
S$.  Write $\partial T = \mathbf p = (p_1, \dotsc, p_n) \subset S$, so that 
$(S, \mathbf p)$ constitutes a punctured surface.  The relative character 
variety $\chi_{G, \mathcal C}(S, \mathbf p)$ carries a natural symplectic 
structure on its smooth stratum.

Let $r : \chi_{G, \mathcal C}(Y, T) \to \chi_{G, \mathcal C}(S, \mathbf p)$ be 
the canonical pullback map (i.e. the map induced by restriction to the 
boundary).  Suppose that $\rho \in \chi_{G, \mathcal C}(Y, T)$ is a regular 
point of $r$.  Then there exists a neighborhood $U \subset \chi_{G, \mathcal 
C}(Y, T)$ of $\rho$ such that the restriction $r\vert_U : U \to \chi_{G, 
\mathcal C}(S, \mathbf p)$ is a Lagrangian embedding.  In particular, if all 
points of of $\chi_{G, \mathcal C}(Y, T)$ are regular, then $r$ is a Lagrangian 
immersion.
\end{thm}

In this paper we restrict our attention to the case of $G = \mathrm{SU}(2)$, 
which we identify in the standard way with the group of unit quaternions:
\begin{displaymath}
\left(\! \begin{array}{rr}
\alpha & \overline \beta \\
-\beta & \overline \alpha
\end{array} \!\right)
= \left(\! \begin{array}{rr}
a + b \mathbf i  & c - d \mathbf i \\
-c - d \mathbf i & a - b \mathbf i
\end{array} \!\right)
\mapsto \alpha + \mathbf j \beta = a + b \mathbf i + c \mathbf j + d \mathbf k,
\end{displaymath}
\begin{equation*} \label{eq:unit-norm}
|\alpha|^2 + |\beta|^2 = a^2 + b^2 + c^2 + d^2 = 1.
\end{equation*}
Here the conjugacy classes are parametrized by the trace 
$\operatorname{tr}(Q)$, which is twice the real part of $Q \in \mathrm{SU}(2)$.
Our choice of conjugacy class is that of traceless elements:
\begin{equation} \label{eq:trless-conj-class}
\mathcal C_0 := \left\{ Q \in \mathrm{SU}(2) \;\middle\vert \; 
\operatorname{tr}(Q) = 0 \right\} = \left\{ \exp \left( \tfrac{\pi}{2} \mathbf 
i \xi \right) \; \middle\vert \; \xi \in \mathfrak{su}(2) \right\}.
\end{equation}
For these choices we adopt the following terminology and notation:
\begin{defn}[Traceless character variety] \label{def:trless-xvar}
Let the pair $(M, T)$ be as above and $\mathcal C_0 \subset \mathrm{SU}(2)$ be 
the conjugacy class of traceless elements. The \textbf{traceless character 
variety} of $M$ relative to $T$ is
\begin{equation} \label{eq:trless-xvar}
\mathscr R(M, T) := \chi_{\mathrm{SU}(2), \mathcal C_0}(M, T).
\end{equation}
\end{defn}

Given a genus-$g$ punctured surface $(\Sigma_g, \mathbf p)$ with $n$ punctures 
$\mathbf p = (p_1, \dotsc, p_n)$, fix meridians $c_1, \dotsc, c_n$ of the 
punctures with $c_j$ oriented according to whether $\varepsilon_j = 
\varepsilon(p_j)$ agrees with the canonical orientation of a point, and form a 
standard presentation
\begin{equation} \label{eq:psurf-std}
\pi_1(\Sigma_g \setminus \mathbf p) = \left\langle a_1, b_1, \dotsc, a_g, b_g, 
c_1, \dotsc, c_n \mid \prod_{i=1}^g [a_i, b_i] \cdot \prod_{j=1}^n 
c_j^{\varepsilon_j} = 1 \right\rangle,
\end{equation}
where $[a_i, b_i]$ denotes the commutator.  Write $\rho(a_i) = A_i$, $\rho(b_i) 
= B_i$ and $\rho(c_j) = C_j$ for each $[\rho] \in \mathscr R(\Sigma_g, \mathbf 
p)$, and collect these into tuples
$\mathbf A, \mathbf B \in G^g$ and $\mathbf C \in \mathcal C_0^n$ in the 
obvious fashion.  In so doing we identify the traceless character variety with
\begin{equation} \label{eq:subvar}
\mathscr R_{g,n} := \left\{ (\mathbf A, \mathbf B, \mathbf C) \in G^g \times 
G^g \times \mathcal C_0^n \mid \prod_{i=1}^g [A_i, B_i] \cdot \prod_{j=1}^n 
C_j^{\varepsilon_j} = 1
\right\} / \text{conj.},
\end{equation}
elements of which we write $[\mathbf A, \mathbf B, \mathbf C]$.

In general, $\mathscr R_{g,n}$ is a stratified symplectic manifold carrying a 
symplectic form $\omega_{g,n}$, called the \emph{Atiyah-Bott-Goldman form}, on 
its top smooth stratum. Under the forgoing identification, the latter is 
precisely the symplectic form of Theorem \ref{thm:lagr-immers}.  When $n$ the 
number of punctures in $\Sigma_g$ is odd, the situation is particularly 
convenient:
\begin{prop}[{\cite[Proposition 2.1]{Hor-sih_trless_xvar}}] 
\label{prop:smooth-punc-parity}
If $n$ is odd, then $(\mathscr R_{g, n}, \omega_{g,n})$ is a smooth symplectic 
manifold of dimension $6g - 6 + 2n$.  If $n$ is even, then it has 
singularities.
\end{prop}

Assume henceforth that $n$ is odd.  The following propositions list those 
symplectic and algebro-topological properties of $\mathscr R_{g, n}$ which are 
needed here:

\begin{prop}[{\cite[Theorem 4.2]{MW-can_bundle_ham_loop_grp_mfld}}] 
\label{prop:surf-monot}
The symplectic manifold $(\mathscr R_{g, n}, \omega_{g,n})$ is 
$\frac{1}{4}$-monotone.
\end{prop}

\begin{prop}[{\cite[Proposition 2.1]{Hor-sih_trless_xvar}}] 
\label{prop:surf-chern}
The minimal Chern number of $\mathscr R_{g, n}$ is equal to $1$.
\end{prop}

\begin{rem} \label{rem:surf-lagr-maslov}
It follows from Proposition \ref{prop:surf-chern} and Equation 
\eqref{eq:maslov2chern} that the minimal Maslov number $N_L = 2$ for any 
compact, simply connected monotone Lagrangian $L \subset (\mathscr R_{g, n}, 
\omega_{g,n})$.
\end{rem}


%% file: sec3-construct.tex
\section{Construction} \label{sec:constr}

Let $Y$ be a connected, closed, oriented $3$-manifold and $K \subset Y$ be an 
$\ell$-component link.  In summary, we will associate to the pair $(Y, K)$ a 
slightly modified balanced multi-pointed Heegaard diagram, the definition of 
which is reviewed at the outset of the construction below.  This diagram 
determines a symplectic manifold, the traceless $\mathrm{SU}(2)$-character 
variety of the Heegaard surface relative to its marked points, as well as a 
pair of Lagrangian submanifolds therein.  The latter correspond to the 
traceless character varieties of the two handlebodies relative to their 
intersection with $K$ and an additional datum $\theta$.  We conclude the 
section by showing that the resulting Lagrangian Floer homology is well defined 
as an abelian group with integer coefficients.

\subsection{Compatible Heegaard diagrams} \label{sec:heeg-diagr}

Recall that a \emph{(balanced) $(2k+m)$-pointed Heegaard diagram} is a tuple 
$\mathcal H = (\Sigma_g, \boldsymbol \alpha, \boldsymbol \beta; \mathbf w, 
\mathbf z, \mathbf q)$ satisfying the following:
\begin{enumerate}[label=(\roman*)] \label{def:bal-ptd-heeg-diagr}

\item $\Sigma_g$ is a closed, oriented surface of genus $g$.

\item $\boldsymbol \alpha = \left( \alpha_1, \dotsc, \alpha_{g + k + m - 1} 
\right)$ and $\boldsymbol \beta = \left( \beta_1, \dotsc, \beta_{g + k + m - 
1} \right)$ are each a collection of disjoint simple closed curves in 
$\Sigma_g$, and each spans a $g$-dimensional lattice in $H_1(\Sigma_g; 
\mathbb Z)$.

\item $\mathbf w = (w_1, \dotsc, w_k)$, $\mathbf z = (z_1, \dotsc, z_k)$ and 
$\mathbf q = (q_1, \dotsc, q_m)$ are mutually disjoint tuples of distinct 
points on $\Sigma_g$.  The points of $\mathbf q$ are called \emph{free 
basepoints}.

\item For $i = 1, \dotsc, k+m$, let $A_i$ denote the connected components of 
$\Sigma_g \setminus \boldsymbol \alpha$ and $B_i$ those of $\Sigma_g \setminus 
\boldsymbol \beta$.  Then for some permutation $\sigma$ of $\left\{ 1, \dotsc, 
k \right\}$, we have $w_i \in A_i \cap B_i$ and $z_i \in A_i \cap
B_{\sigma(i)}$ for $i = 1, \dotsc, k$, and $q_j \in A_j \cap B_j$ for $j = 1, 
\dotsc, m$.

\end{enumerate}
(This definition can be found e.g. in \cite[Section 
3.1]{MO-heeg_floer_hom_int_surg_link} with different notation around free 
  basepoints.)  We will construct such a diagram for $(Y, K)$, with $k \ge 
\ell$ the relative bridge number of $K$ with respect to $\Sigma_g$.  An initial 
choice of a single free basepoint will be thickened into a triple, in order to 
arrive at an odd number of punctures.

In the following let $\mathrm{Crit_i}(f)$ denote the set of critical points of 
a Morse function $f$ having Morse index $i \in \mathbb Z_{\ge 0}$.  To 
construct the desired Heegaard diagram, begin by fixing a Morse function $f_K : 
K \to \mathbb R$ with $\mathrm{Crit}_0(f_K) = \left\{ a_1, \dotsc, a_k 
\right\}$ and $\mathrm{Crit}_1(f_K) = \left\{ b_1, \dotsc, b_k \right\}$ for $k 
\geq \ell$.  Extend $f_K$ to a self-indexing Morse function $f_{(Y, K)} : Y \to 
\mathbb R$ so that, for some pair of points $a_0, b_0 \in Y \setminus K$,
\begin{enumerate}[label=(\roman*)]
\item $\mathrm{Crit}_0(f_{(Y, K)}) = \mathrm{Crit}_0(f_K) \cup \left\{ a_0 
\right\}$;
\item $\mathrm{Crit}_3(f_{(Y, K)}) = \mathrm{Crit}_1(f_K) \cup \left\{ b_0 
\right\}$.
\end{enumerate}
The function $f_{(Y, K)}$ determines the components of a balanced 
$(2k+1)$-pointed Heegaard diagram in the standard fashion:  The Heegaard 
surface is the level set $\Sigma_g = f_{(Y, K)}^{-1}\left( \frac 3 2 \right)$.  
The stable manifolds of the $g+k$ index-$1$ critical points intersect 
$\Sigma_g$ in a $(g+k)$-tuple of disjoint simple closed curves $\boldsymbol 
\alpha = \left\{ \alpha_1, \dotsc, \alpha_{g+k} \right\}$, the 
$\alpha$-attaching circles, which determine a handlebody $U_\alpha$ with 
boundary $\Sigma_g$.  Similarly the unstable manifolds of the index-$2$ 
critical points intersect $\Sigma_g$ in the $\beta$-attaching circles 
$\boldsymbol \beta = \left\{ \beta_1, \dotsc, \beta_{g+k} \right\}$, 
determining the handlebody $U_\beta$.  The intersection $\Sigma_g \cap K$ is a 
set of $2k$ points.  By choosing an orientation of $K$, some ordering of its 
components and an arbitrary starting point on each, we obtain an indexing $w_1, 
z_1, \dotsc, w_k, z_k$ of the points in $\Sigma_g \cap K$ as described above.  
Any choice of flowline connecting the critical points disjoint from $K$, viz. 
$a_0$ and $b_0$, intersects $\Sigma_g$ in the single free basepoint.

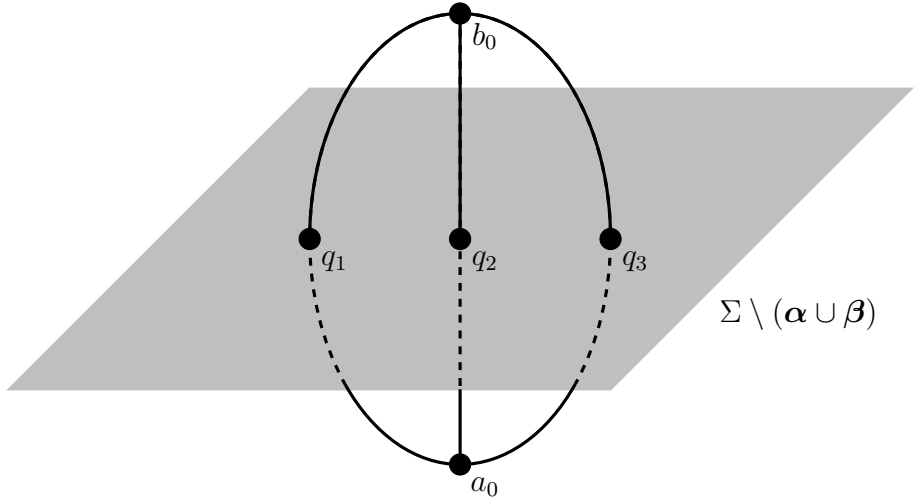
\begin{figure} \label{fig:theta-graph}
\begin{tikzpicture}
\filldraw[lightgray, thick] (0, 2) -- (4, 6) -- (12, 6) -- (8, 2) -- cycle;

\coordinate (0) at (6, 7);
\coordinate (1) at (4, 4);
\coordinate (2) at (6, 4);
\coordinate (3) at (8, 4);
\coordinate (4) at (6, 1);

\node[] at (10.5, 3) { $\Sigma \setminus (\boldsymbol \alpha \cup \boldsymbol 
\beta)$ };
\path (0) +(0.33, -0.3) node[] (b) { $b_0$ };
\path (1) +(0.33, -0.3) node[] (x1) { $q_1$ };
\path (2) +(0.33, -0.3) node[] (x2) { $q_2$ };
\path (3) +(0.33, -0.3) node[] (x3) { $q_3$ };
\path (4) +(0.33, -0.3) node[] (a) { $a_0$ };

\filldraw[black] (0) circle (4pt);
\filldraw[black] (1) circle (4pt);
\filldraw[black] (2) circle (4pt);
\filldraw[black] (3) circle (4pt);
\filldraw[black] (4) circle (4pt);

\draw[black, very thick] (0) -- (2);
\draw[black, dashed, very thick] (0) -- (6, 2);
\draw[black, very thick] (6, 2) -- (4);

\draw[black, very thick] (0) arc (90:180:2 and 3);
\draw[black, very thick] (0) arc (90:0:2 and 3);
\draw[black, very thick, dashed] (0) arc (90:220:2 and 3);
\draw[black, very thick, dashed] (0) arc (90:-40:2 and 3);
\draw[black, very thick] (4) arc (270:220:2 and 3);
\draw[black, very thick] (4) arc (270:320:2 and 3);
\end{tikzpicture}
\caption{Theta graph $\theta$ in a neighborhood of $Y \setminus K$.}
\end{figure}

To suit our purposes, we modify this construction at the last step.  
Specifically, we choose three distinct flowlines connecting $a_0$ and $b_0$ 
rather than one, and collect their intersections with $\Sigma_g$ into the 
triple $\mathbf q = (q_1, q_2, q_3)$.  With regard to the Heegaard diagram 
above, one may think of this as thickening the single free basepoint into a 
triple, the auxiliary datum anticipated above.  We will work with the 
$(2k+3)$-pointed Heegaard diagram $\mathcal H = (\Sigma_g, \boldsymbol \alpha, 
\boldsymbol \beta; \mathbf w, \mathbf z, \mathbf q)$ and refer to any
diagram produced in this manner as being compatible with $(Y, K)$.  Denote the 
union of the three flowlines chosen above by $\theta \subset Y \setminus K$, 
referred to as a theta graph; see fig. \ref{fig:theta-graph} for a schematic 
illustration.

Note that this process of deriving a compatible $(2k+3)$-pointed Heegaard 
diagram from a link is reversible.  We may assume that the indexing of the 
connected components in (iii) of the definition \ref{def:bal-ptd-heeg-diagr} 
agrees with the one chosen for the points $\mathbf w, \mathbf z$.  For $i = 1, 
\dotsc, k$, draw a simple curve in the connected component $A_i \subset 
\Sigma_g$ connecting $w_i$ and $z_i$, and then push the interior of this arc 
into $U_\alpha$.  Likewise connect $z_{\sigma^{-1}(i)}$ to $w_i$ by an arc in 
$B_i$ and push this off into $U_\beta$.  The union of all the pushed-off arcs 
forms a link in $Y = U_\alpha \cup_{\Sigma_g} U_\beta$.  Now choose a point in 
each of $A_{k+1}$ and $B_{k+1}$, and connect the points of $\mathbf q$ to each 
of these by three simple arcs, disjoint except at the common endpoint.  Pushing 
the triple of arcs in $A_{k+1}$ off into $U_\alpha$ and that in $B_{k+1}$ off 
into $U_\beta$ produces a theta graph $\theta \in Y$.

For brevity, fix the notation $K_\square = U_\square \cap K$ and 
$\theta_\square = U_\square \cap \theta$ for $\square = \alpha$ or $\beta$.

\begin{rem}[] \label{rem:triv-hbdy-tangle}
An important aspect of such a Heegaard splitting for our purposes is that the 
pairs $(U_\square, K_\square \cup \theta_\square)$ are trivial:  $K_\square$ is 
consists of properly embedded arcs which are boundary parallel.  The tripod 
$\theta_\square$ is similarly boundary parallel in the sense that there is a 
simultaneous isotopy of the closure of each of its three strands onto 
$\Sigma_g$.
\end{rem}

\subsection{Lagrangians from Heegaard diagrams} \label{sec:lagr-heeg-diagr}

Since $\# \mathbf w \cup \mathbf z \cup \mathbf q = 2k + 3$ is odd, the 
traceless $\mathrm{SU}(2)$-character variety $\mathscr R_{g, 2k+3} = \mathscr 
R(\Sigma_g, \mathbf w \cup \mathbf z \cup \mathbf q)$ is a smooth symplectic 
manifold by Proposition \ref{prop:smooth-punc-parity}.  The attaching circles 
determine subvarieties in $\mathscr R_{g, 2k+3}$ in a natural way:
\begin{equation} \label{eq:lagr-def}
L_\alpha := \left\{ [\rho] \in \mathscr R_{g, 2k+3} \mid \rho([\alpha_i]) = 1, 
\text{ for } i = 1, \dotsc, g + k \right\}
\end{equation}
and
\begin{equation} \label{eq:lagr-def}
L_\beta := \left\{ [\rho] \in \mathscr R_{g, 2k+3} \mid \rho([\beta_i]) = 1, 
\text{ for } i = 1, \dotsc, g + k \right\}.
\end{equation}
Let $c_i, d_i$ be meridians of the punctures $w_i, z_i$, resp., and $m_1, m_2, 
m_3$ those of $q_1, q_2, q_3$, resp.  For simplicity we use the same notation 
for their homotopy classes relative to some choice of basepoint.  Fixing an 
orientation of $K$ determines an orientation of the points $w_i, z_i$, and this 
in turn orients their meridians.  We can also fix a coherent orientation of 
$m_1, m_2, m_3$, so that we obtain the standard presentation of $\pi_1(\Sigma_g 
\setminus (\mathbf w \cup \mathbf z \cup \mathbf q))$ :
\begin{equation} \label{eq:std-pres}
\left\langle \begin{array}{c} a_1, b_1, \dotsc, a_g, b_g, \\ c_1, d_1, \dotsc 
c_k, d_k, \\ m_1, m_2, m_3 \end{array} \;\middle\vert\; \left( \prod_{i=1}^g 
[a_i, b_i] \right) \cdot \left( \prod_{j=1}^k c_j \cdot d_j^{-1} \right) \cdot 
(m_1 \cdot m_2 \cdot m_3) = 1 \right\rangle.
\end{equation}
It is evident that, if a representation $\rho$ sends all of the $[\alpha_i]$ to 
$1 \in \mathrm{SU}(2)$, then $\rho(c_i) = -\rho(d_i)$ and $\rho(m_1) \rho(m_2) 
\rho(m_3) = 1$.

It will be useful to identify these subvarieties with the traceless character 
varieties of the pairs $(U_\alpha, K_\alpha \cup \theta_\alpha)$ and $(U_\beta, 
K_\beta \cup \theta_\beta)$, resp.  We make this precise for $L_\alpha$ in the 
following; the argument for $L_\beta$ is essentially the same.

\begin{lem}[] \label{lem:hbdy-xvar}
Let $\widetilde U_\alpha \subset U_\alpha$ be the $3$-manifold with boundary 
$\partial \widetilde U_\alpha = \Sigma_g \sqcup S$ where $S$ is the $2$-sphere 
boundary of a small $3$-ball neighborhood of the critical point $a_0$ (the 
vertex of the theta graph $\theta$ lying in $U_\alpha$), and write $\widetilde 
\theta_\alpha = \widetilde U_\alpha \cap \theta$.  Let
\begin{displaymath}
r : \mathscr R\left( \widetilde U_\alpha, K_\alpha \cup \widetilde 
\theta_\alpha \right) \to \mathscr R(\Sigma_g, \mathbf w \cup \mathbf z \cup 
\mathbf q) \times \mathscr R(S, S \cap \theta),
\end{displaymath}
be the pullback map.  Then $r$ is injective, and its image in the first factor 
of the target is naturally identified with $L_\alpha$.
\end{lem}
\begin{proof}
A computation shows that $\mathscr R(S, S \cap \theta) = \mathscr R_{0,3}$ is 
just the point $\left\{ [\mathbf i, \mathbf j, -\mathbf k] \right\}$ 
\cite[Example 2.3]{Hor-sih_trless_xvar}, so the target of $r$ is naturally 
identified with $\mathscr R_{g,2k+3}$.  It is easy to see that $L_\alpha$ is 
composed of classes of precisely those traceless representations on $(\Sigma_g, 
\mathbf w \cup \mathbf z \cup \mathbf q)$ which extend over $(\widetilde 
U_\alpha, K_\alpha \cup \widetilde \theta_\alpha)$.  Because the canonical map 
$\pi_1(\Sigma_g \setminus (\mathbf w \cup \mathbf z \cup \mathbf q)) \to 
\pi_1(\widetilde U_\alpha \setminus (K_\alpha \cup \widetilde \theta_\alpha))$ 
is surjective, injectivity follows ultimately from the left-exactness of the 
contravariant functor $\Hom(-, \mathrm{SU}(2))$.
\end{proof}

This characterization and Remark \ref{rem:triv-hbdy-tangle} lead to the 
conclusion that, for a fixed genus and relative bridge number, the subvarieties 
arising from a compatible Heegaard diagram as above are equivalent in the 
following sense.

\begin{lem}[] \label{lem:diffeom-equiv}
Let $\mathcal H$ be an arbitrary balanced $(2k+3)$-pointed Heegaard diagram and 
$L_\alpha, L_\beta \subset \mathscr R_{g,2k+3}$ the subvarieties induced as 
above.  Then there exists a symplectomorphism $\phi : \mathscr R_{g,2k+3} \to 
\mathscr R_{g,2k+3}$ such that $\phi(L_\alpha) = L_\beta$.
\end{lem}
\begin{proof}
By Lemma \ref{lem:hbdy-xvar} each of $L_\alpha$ and $L_\beta$ is identified 
with the images of$\mathscr R(U_\alpha, K_\alpha \cup \theta_\alpha)$ and 
$\mathscr R(U_\beta, K_\beta \cup \theta_\beta)$, resp., in $\mathscr 
R_{g,2k+3}$.  As described in Remark \ref{rem:triv-hbdy-tangle}, the pairs 
$(U_\alpha, K_\alpha \cup \theta_\alpha)$ and $(U_\beta, K_\beta \cup 
\theta_\beta)$ are trivial, and so there exists a diffeomorphism $(U_\beta, 
K_\beta \cup \theta_\beta) \to (U_\alpha, K_\alpha \cup \theta_\alpha)$ 
restricting to a pointed diffeomorphism $\phi$ of their common boundary 
$\Sigma_g$.  Since mapping classes act by symplectomorphisms, this in turn 
induces a symplectomorphism $\phi : \mathscr R_{g,2k+3} \to \mathscr 
R_{g,2k+3}$ which maps the image of $\mathscr R(U_\alpha, K_\alpha \cup 
\theta_\alpha)$ to that of $\mathscr R(U_\beta, K_\beta \cup \theta_\beta)$, 
i.e. $\phi(L_\alpha) = L_\beta$.
\end{proof}

\begin{prop}[] \label{prop:lagr-embed}
The subvarieties $L_\alpha, L_\beta$ are embedded Lagrangian submanifolds in 
$(\mathscr R_{g,n}, \omega_{g,n})$.
\end{prop}
\begin{proof}
By Theorem \ref{thm:lagr-immers} the boundary restriction map $r$ of Lemma 
\ref{lem:hbdy-xvar} is a Lagrangian immersion at its regular points, which are 
precisely those mapped to the (classes of) irreducible representations in 
$\mathscr R_{g,2k+3}$.  Since there are no reducibles in the latter, $r$ is a 
global immersion.  By Lemma \ref{lem:hbdy-xvar} it is injective, hence its 
image, identified with $L_\alpha$, is an embedded Lagrangian.
\end{proof}

We conclude the subsection by identifying the diffeomorphism type of $L_\alpha$ 
and $L_\beta$.

\begin{lem} \label{lem:sph-prod}
The subvarieties $L_\alpha, L_\beta$ are diffeomorphic to the product of 
spheres $\left( S^3 \right)^g \times \left( S^2 \right)^k$.
\end{lem}
\begin{proof}
By Lemma \ref{lem:diffeom-equiv} it suffices to prove the result for $L_\alpha$ 
where the $\alpha$-curves are standard in the following sense:  The curves 
$\alpha_1, \dotsc, \alpha_g$ represent the generators $a_1, \dotsc, a_g$ 
appearing in the standard presentation of $\mathscr R_{g,2k+3}$ described in 
Section \ref{sec:xvar-surf}.  Let $c_j, d_j$ for $j = 1, \dotsc, k$ and $m_1, 
m_2, m_3$ be as above and collect the images $\rho(c_j)$, $\rho(d_j)$ and 
$\rho(m_p)$ into tuples $\mathbf B \in G^g$, $\mathbf C \in \mathcal C_0^k$ and 
$\mathbf M \in \mathcal C_0^3$, resp., where $G = \mathrm{SU}(2)$.  Then 
elements of $L_\alpha$ may be expressed as $[\mathbf 1, \mathbf B, \mathbf C, 
\mathbf M] \in G^g \times G^g \times \mathcal C_0^k \times \mathcal C_0^3 / G$ 
subject to the constraint $M_1 M_2 M_3 = 1$.

Consider the map $G^g \times \mathcal C_0^k \to L$ defined by
\begin{equation} \label{eq:lagr-proj}
(\mathbf B, \mathbf C) \mapsto [\mathbf 1, \mathbf B, \mathbf C, (\mathbf i, 
\mathbf j, -\mathbf k)].
\end{equation}
This map is injective since the common stabilizer of $\mathbf i, \mathbf j, 
-\mathbf k$ is the center $\left\{ \pm 1 \right\}$ of $\mathrm{SU}(2)$.  On the 
other hand, for any $[\mathbf 1, \mathbf B, \mathbf C, \mathbf M] \in L$, there 
is a unique $R \in \mathrm{SO}(3)$ such that $R \mathbf M R^{-1} = (\mathbf i, 
\mathbf j, -\mathbf k)$, hence
\begin{equation} \label{eq:surj}
[\mathbf 1, \mathbf B, \mathbf C, \mathbf M] = [\mathbf 1, R \mathbf B R^{-1}, 
R \mathbf C R^{-1}, (\mathbf i, \mathbf j, -\mathbf k)],
\end{equation}
is the image under the map \eqref{eq:lagr-proj} of $R(\mathbf B, \mathbf 
C)R^{-1}$, demonstrating that this map is surjective.  Moreover, $R$ varies 
smoothly with $\mathbf M$, and therefore \eqref{eq:lagr-proj} is a 
diffeomorphism.
\end{proof}

An easy consequence of Lemma \ref{lem:sph-prod} is the following:

\begin{cor}
\label{cor:spin-lagr}
The Lagrangian submanifolds $L_\alpha, L_\beta$ are spin.
\end{cor}
\begin{proof}
All spheres are spin manifolds, and the product of spin manifolds is again 
spin, as a consequence of the Whitney sum formula.
\end{proof}

\subsection{Floer homology} \label{sec:sik-hf}
In general the Lagrangian intersection $L_\alpha \cap L_\beta$ will not be 
transverse.  However, standard arguments show that a Hamiltonian isotopy may be 
chosen generically after which transversality is achieved, and that such 
isotopies produce canonically isomorphic Lagrangian Floer homologies.  It may 
in fact be possible to accomplish this by means of the \emph{holonomy 
perturbations} defined in \cite[Section 6.2]{CHK-tngl_xvar}, which are 
analogous to the perturbations employed by Kronheimer and Mrowka in 
\cite{KM-knot_hom_grp_inst} for singular instanton knot homology, although 
further study into this is needed.  In any event we assume henceforth that 
$L_\alpha$ and $L_\beta$ intersect transversely.

We now proceed to prove Theorem \ref{thm:A} of the introduction, recapitulated 
in the following:
\begin{thm}
\label{thm:sik-hf}
The Lagrangian Floer homology $\mathrm{HF}(L_\alpha, L_\beta)$ is well-defined.  
Moreover, the pair $(L_\alpha, L_\beta)$ is relatively spin, so that the 
homology may be taken with coefficients in $\mathbb Z$.
\end{thm}

The first step is to establish that we are in the monotone setting, so that 
Oh's results for this case summarized in the previous section can be applied.

\begin{prop}
\label{prop:monot-lagr}
The Lagrangians $L_\alpha, L_\beta$ are monotone in $(\mathscr R_{g,2k+3}, 
\omega_{g,2k+3})$ and have minimal Maslov number $N_L = 2$.
\end{prop}
\begin{proof}
By Lemma \ref{lem:sph-prod} both $L_\alpha, L_\beta$ are compact and simply 
connected, and therefore monotone by Proposition \ref{prop:sc-monot}.  As for 
the second claim, recall from Proposition \ref{prop:surf-chern} that the 
minimal Chern number of $(\mathscr R_{g,n}, \omega_{g,n})$ is $N_M = 1$ for any 
$n$, hence $N_L = 2$ by Equation \eqref{eq:maslov2chern}.
\end{proof}

\begin{proof}[Proof of Theorem \ref{thm:sik-hf}]
By Lemma \ref{lem:sph-prod} and Proposition \ref{prop:monot-lagr}, the pair 
$(L_\alpha, L_\beta)$ satisfy the conditions of Theorem \ref{thm:hf-gw-invar}.  
Since, as demonstrated in Lemma \ref{lem:diffeom-equiv}, there is a 
symplectomorphism of $\mathscr R_{g,2k+3}$ taking $L_\alpha$ to $L_\beta$, we 
can apply Corollary \ref{cor:hf-gw-sympl} to conclude that the Lagrangian Floer 
homology $\mathrm{HF}(L_\alpha, L_\beta)$ is well defined.  Being that 
$L_\alpha$ and $L_\beta$ are both spin by Corollary \ref{cor:spin-lagr}, they 
are trivially relatively spin in $\mathscr R_{g,2k+3}$ as a pair, hence their 
Floer homology may be taken with coefficients in $\mathbb Z$.
\end{proof}

\subsubsection{Gradings} \label{sec:grd}
For any compact monotone Lagrangians $L_0, L_1$ in a symplectic manifold $M$ 
with $N_L \ge 2$, the Lagrangian Floer homology $\mathrm{HF}(L_0, L_1)$ is 
endowed with a relative $\mathbb Z/N_L$-grading.  When $M, L_0, L_1$ are 
orientable, its reduction $\!\!\!\!\mod 2$ is equivalent to the grading which 
results from fixing orientations of the Lagrangians and then assigning a degree 
in $\left\{ \pm 1 \right\}$ to each point $x \in L_0 \cap L_1$ according to 
whether the induced orientation agrees with the canonical orientation of the 
point $x \in M$.  Thus $\mathrm{HF}(L_0, L_1)$ categorifies the algebraic 
intersection number of $L_0, L_1 \subset M$.

Since $N_L = 2$ for our Lagrangians $L_\alpha, L_\beta$, and these are 
certainly orientable, the standard $\mathbb Z/N_L$-grading on 
$\mathrm{SIK}(\mathcal H) = \mathrm{HF}(L_\alpha, L_\beta)$ obtained after 
fixing a base intersection point in grading $0$ is exactly the $\mathbb 
Z/2$-grading induced by fixing orientations, up to an overall sign in the 
degrees.  Finer $\mathbb Z/N$-gradings, where $N$ is divisible by $N_L = 2$, 
are possible \textit{a priori}.  Defining such a grading requires an 
understanding of the space $\pi_2(L_\alpha, L_\beta)$ of topological Whitney 
disks with Lagrangian boundary conditions connecting intersection points, and 
we plan to visit this subject in future work on $\mathrm{SIK}$.

%% file: sec4-top-invar.tex
\section{Topological invariance} \label{sec:top-invar}
In this section we prove Theorem \ref{thm:B} of the introduction, viz. that the 
isomorphism class of the homology $\sik(\mathcal H)$ constructed above is a 
topological invariant of the pair $(Y, K)$.  To this end we make use of the 
fact that any two compatible diagrams $\mathcal H$ and $\mathcal H'$ are 
related by a sequence of multi-pointed Heegaard moves, as described by the 
following.
\begin{prop}[{\cite[Proposition 
3.9]{MO-heeg_floer_hom_int_surg_link}}] \label{prop:heeg-mvs}
Let $\mathcal H$ and $\mathcal H'$ be two multi-pointed Heegaard diagrams 
compatible with the pair $(Y, K)$ (in the sense of Section \ref{sec:heeg-diagr} 
above).  Then $\mathcal H$ can be transformed into $\mathcal H'$ by a sequence 
of the following moves.
\begin{enumerate}[label=(\roman*)]
\item \label{def:isotopy} An \textbf{isotopy} of an attaching circle supported 
in the complement of the punctures $\mathbf w \cup \mathbf z \cup \mathbf q$;

\item \label{def:handleslide} A \textbf{handleslide} among the $\alpha$- or 
among the $\beta$-attaching circles which is supported in the complement of 
$\mathbf w \cup \mathbf z \cup \mathbf q$;

\item \label{def:stab-surf} \textbf{surface (de-)stabilization}:  The new 
Heegaard diagram
\begin{equation} \label{eq:stab-heeg-diagr}
\mathcal H' = (\Sigma_{(g + 1)}, \boldsymbol \alpha \cup \alpha', \boldsymbol 
\beta \cup \beta'; \mathbf w, \mathbf z, \mathbf q)
\end{equation}
is formed by taking the connected sum of the Heegaard surface $\Sigma_g$ in 
$\mathcal H$ with a torus along a disk which is disjoint from the $\alpha$- and 
$\beta$-attaching circles and from the marked points of $\mathbf w \cup \mathbf 
z \cup \mathbf q$.  On the torus are two new attaching curves $\alpha'$ and 
$\beta'$ which intersect transversely in a single point.  (Resp. the inverse of 
this procedure, with the roles of $\mathcal H$ and $\mathcal H'$ exchanged.)

\item \label{def:stab-link} \textbf{link (de-)stabilizations}:  The new 
Heegaard diagram
\begin{equation} \label{eq:stab-heeg-diagr}
\mathcal H' = (\Sigma_g, \boldsymbol \alpha \cup \alpha', \boldsymbol \beta 
\cup \beta'; \mathbf w \cup \left\{ w' \right\}, \mathbf z \cup \left\{ z' 
\right\}, \mathbf q)
\end{equation}
is formed by the addition to $\mathcal H$ of a new pair of attaching circles 
$\alpha'$ and $\beta'$ and of two new marked points $w'$ and $z'$, satisfying 
the following:  The circles $\alpha', \beta'$ are disjoint from the original 
ones of $\boldsymbol \alpha, \boldsymbol \beta$, resp., and intersect 
transversely in two points.  They bound disks $A', B' \subset \Sigma_g$, resp., 
such that $w' \in A' \cap B'$ while $z' \in A' \setminus B'$ and there is 
precisely one of the original marked points $z_0 \in \mathbf z$ which lies in 
$B' \setminus A'$.  (Resp. the inverse of this procedure.)
\end{enumerate}
\end{prop}

\begin{rem} \label{rem:heeg-mv-diff}
What we have called \emph{surface stabilization} in (iii) above is referred to 
as \emph{index-one/two stabilization} in \cite[Proposition 
3.9]{MO-heeg_floer_hom_int_surg_link}.  Moreover, (iv) link stabilization is 
  one of two kinds of \emph{index-zero/three stabilization} enumerated there by 
the authors, along with \emph{free stabilization}.  The latter move serves to 
add a new free basepoint, of the kind which we have thickened into the triple 
$\mathbf q$.  A well-defined Lagrangian Floer theory could result from certain 
free stabilizations, so long as the total number of marked points remains odd.  
For the moment we have nothing to motivate such moves, nor are they needed to 
obtain the equivalence classes of compatible Heegaard diagrams we have used in 
this construction, and we therefore exclude free stabilization from our
consideration.

Note also that we have dispensed with the requirement that the Heegaard moves 
in question be \emph{admissible}, in the sense given there by the authors, as 
this is not relevant in our context, and the weaker statement is clearly no 
less valid.
\end{rem}

In light of Proposition \ref{prop:heeg-mvs}, it will suffice to prove the 
following:
\begin{thm} \label{thm:heeg-mv-invar}
Let $\mathcal H$ and $\mathcal H'$ be two multi-pointed Heegaard diagrams 
compatible with the pair $(Y, K)$.
\begin{enumerate}
  \item If $\mathcal H$ and $\mathcal H'$ are related by an isotopy or 
    handleslide, then the respective Lagrangian pairs $(L_\alpha, L_\beta)$ and 
    $(L_\alpha', L_\beta')$ determined by these diagrams (as in Section 
    \ref{sec:lagr-heeg-diagr}) are identical, hence $\sik(\mathcal H) = 
    \sik(\mathcal H')$.
  \item If $\mathcal H$ and $\mathcal H'$ are related by surface or link 
stabilization, then $\sik(\mathcal H)$ and $\sik(\mathcal H')$ are canonically 
isomorphic as relatively $\mathbb Z/2$-graded abelian groups.
\end{enumerate}
\end{thm}

The following subsection is devoted to proving the forgoing theorem.  With this 
result in place, we will define $\mathrm{SIK}(Y, K)$ to be the isomorphism 
class of the Lagrangian Floer homology $\mathrm{SIK}(\mathcal H)$ for any 
choice of compatible diagram $\mathcal H$.

\subsection{Invariance under multi-pointed Heegaard moves} 
\label{sec:invar-heeg-mv}
The next proposition proves statement (1) of Theorem \ref{thm:heeg-mv-invar}.
\begin{prop} \label{prop:isot-hndlsld-invar}
Let the Heegaard diagrams $\mathcal H$ and $\mathcal H'$ be related by an 
isotopy \ref{def:isotopy} or handleslide \ref{def:handleslide}, and let 
$(L_\alpha, L_\beta)$ and $(L_\alpha', L_\beta')$ be their respective induced 
Lagrangian pairs, as defined in Section \ref{sec:constr}.  Then the latter are 
identical Lagrangian pairs in $\mathscr R_{g,n}$, i.e. $L_\alpha = L_\alpha'$ 
and $L_\beta = L_\beta'$, and consequently $\mathrm{SIK}(\mathcal H) = 
\mathrm{SIK}(\mathcal H')$.
\end{prop}
\begin{proof}
The case of an isotopy is trivial, as the image of a loop under an element of 
the character variety is determined by the former's free homotopy class.

For a handleslide, we may assume without loss of generality that the curve 
$\alpha_1$ of $\mathcal H$ is slid over $\alpha_2$ to produce $\alpha_1'$ in 
$\mathcal H'$.  By hypothesis these three curves bound a pair of pants in 
$\Sigma_g \setminus (\mathbf w \cup \mathbf z \cup \mathbf q)$, hence any 
$\mathrm{SU}(2)$-representation taking two of them to $1$ does likewise for the 
third.  It follows that $L_\alpha = L_\alpha'$.
\end{proof}

In the remainder of this section we verify statement (2) of Theorem 
\ref{thm:heeg-mv-invar}.  We will establish the claimed the canonical 
isomorphism by associating to the Heegaard move in question a kind of Cerf 
decomposition of the triple $(Y, K, \theta)$, then deriving from this 
decomposition a generalized Lagrangian correspondence among traceless character 
varieties of punctured surfaces and finally carrying out computations on the 
quilted Floer homology \cite{WW-qltd_floer_cohom} of the correspondence.  This 
largely follows the work of Wehrheim and Woodward \cite{WW-floer_fld_thry_tngl} 
on quilted-Floer invariants of tangles between punctured surfaces, as made 
precise in the following:
\begin{defn}[] \label{def:tangle}
Let $(\Sigma_-, \mathbf p_-)$ and $(\Sigma_+, \mathbf p_+)$ be punctured 
surfaces, and let the decorations $\overline{\Sigma_-}$ and $\overline{\mathbf 
p_-}$ indicate reversal of orientation.  A tangle from $(\Sigma_-, \mathbf 
p_-)$ to $(\Sigma_+, \mathbf p_+)$ is a tangle $(X, T)$ together with an 
orientation-preserving diffeomorphism $\phi : \partial X \to 
\overline{\Sigma_-} \sqcup \Sigma_+$ which restricts to an 
orientation-preserving bijection $\phi\vert_{\partial T} : \partial T \to 
\overline{\mathbf p} \sqcup \mathbf p'$.
A tangle $(X, T)$ is called
\begin{itemize}
\item \textbf{elementary} if there is an orientation-preserving diffeomorphism 
$X \to \Sigma_g \times [0, 1]$, where $\Sigma_g$ is an orientable genus-$g$ 
surface;
\item \textbf{cylindrical} if there is an orientation-preserving diffeomorphism 
$(X, T) \to (\Sigma_g, \mathbf p) \times [0, 1]$, where the $(\Sigma_g, \mathbf 
p)$ is a punctured surface.
\end{itemize}
\end{defn}

In the following we use a Morse function to obtain a decomposition of $(Y, K 
\cup \theta)$ into tangles $(Y_{i(i+1)}, Y_{i(i+1)} \cap (K \cup \theta))$ with 
oriented boundary surfaces $\partial Y_{i(i+1)} = \overline{\Sigma_i} \sqcup 
\Sigma_{i+1}$ punctured where they intersect $K \cup \theta$.  It will be 
convenient to suppress the intersections with $K \cup \theta$ from the notation 
by writing $Y_{i(i+1)}^*$ for the tangle and e.g. $\Sigma_i^*$ for the 
punctured surface.

Let $\mathcal H$ be a compatible Heegaard diagram induced by a Morse function 
$f$ and $\mathcal H'$ the result of surface or link stabilization on $\mathcal 
H$.  Take $a_\pm = f(\Sigma_g) \pm \varepsilon$ for some $\varepsilon > 0$.  
When $\varepsilon$ is taken sufficiently small, the preimage $N = f^{-1}([a_-, 
a_+])$ is a normal neighborhood of the Heegaard surface $\Sigma_g$ in $\mathcal 
H$ which contains no critical points of $f$.  The normalized gradient flow of 
$f$ serves to identify the punctured surfaces $\Sigma_1^* = f^{-1}(a_-) \cap (K 
\cup \theta)$ and $\Sigma_3^* = f^{-1}(a_+) \cap (K \cup \theta)$ each with 
$(\Sigma_g, \mathbf w \cup \mathbf z \cup \mathbf q)$, so that we may consider 
their traceless character varieties to be identified with each other as 
$\mathscr R_{g,2k+3}$.

We may deform the restriction of $f$ to the interior of $N$ so as to produce a 
new Morse function $f'$ which induces the stabilized Heegaard diagram $\mathcal 
H'$.  We write its punctured Heegaard surface as $\Sigma_2^*$.  Then the 
function $f'$ induces the tangle decomposition
\begin{equation} \label{eq:tangle-decomp}
(Y, K \cup \theta) = B_- \cup_{\Sigma_0^*} Y_{01}^* \cup_{\Sigma_1^*} Y_{12}^* 
\cup_{\Sigma_2^*} Y_{23}^* \cup_{\Sigma_3^*} Y_{34}^* \cup_{\Sigma_4^*} B_+,
\end{equation}
where each of $\Sigma_0^*$ and $\Sigma_4^*$ is a $2$-sphere, disjoint from $K$ 
and meeting $\theta$ in three punctures, which bounds the small $3$-ball 
neighborhood $B_-$ or $B_+$, resp., of a vertex of $\theta$.

The tangles of \eqref{eq:tangle-decomp} bear the following descriptions.  The 
``outer pieces'' $B_- \cup Y_{01}^*$ and $Y_{34}^* \cup B_+$ can be identified 
in an obvious way with the pairs $(U_\alpha, K_\alpha \cup \theta_\alpha)$ and 
$(U_\beta, K_\beta \cup \theta_\beta)$, resp.  Each of $Y_{01}, Y_{34}$ meets 
$\theta$ in three trivial arcs running from one boundary component to the 
other, and meets $K$ in $k$ trivial arcs with both endpoints on the same 
boundary component: on the outgoing boundary $\partial^+Y_{01} = \Sigma_1$, 
where we refer to the arcs as \emph{cups}, and on the incoming boundary 
$\partial^-Y_{34} = \overline{\Sigma_3}$, where we refer to them as 
\emph{caps}.

The ``inner piece'' $Y_{12}^* \cup Y_{23}^*$ is diffeomorphic to the 
cylindrical tangle $(\Sigma_g, \mathbf w \cup \mathbf z \cup \mathbf q) \times 
[0, 1]$, containing $2k+3$ trivial arcs running from $\Sigma_1$ to $\Sigma_3$.  
The characterizations of the individual constituent tangles $Y_{12}^*$ and 
$Y_{23}^*$ depend on the type of stabilization performed.  In the case of 
surface stabilization, they can be viewed as the result of taking the 
cylindrical tangle $\Sigma_2^* \times [0, 1]$ and attaching a $2$-handle along 
$\alpha'$ or $\beta'$, resp.  In the case of link stabilization, they are both 
elementary and each contain $2k+3$ arcs running between the two boundary 
components.  On the one hand, $Y_{12}^*$ has in addition a cup connecting $w', 
z' \in \Sigma_2$, and one of the other arcs connects $z_0 \in \Sigma_1$ with 
$z_0 \in \Sigma_2$.  On the other $Y_{23}^*$ has a cap connecting $z_0, w' \in 
\Sigma_2$ and $z' \in \Sigma_2$ is connected by an arc to $z_0 \in \Sigma_3$.

For $i = 0, \dotsc, 4$, let $M_i = \mathscr R(\Sigma_i^*)$ denote the traceless 
character variety of the punctured surface, and $M_i^-$ the result of replacing 
the Atiyah-Bott-Goldman symplectic form by its negative.  Recall that, for $i = 
1, \dotsc, 4$, the traceless character variety $\mathscr R(\partial 
Y_{(i-1)i}^*)$ of the disconnected boundary, comprising the pair of punctured 
surfaces $\Sigma_{i-1}^*$ and $\Sigma_i^*$, is by definition the product 
$M_{i-1} \times M_i$.  Write $r_{(i-1)i} : \mathscr R(Y_{(i-1)i}^*) \to 
\mathscr R(\partial Y_{(i-1)i}^*)$ for the pullback (or boundary-restriction) 
map.  Then $L_{(i-1)i} := \operatorname{im}(r_{(i-1)i}) \subset M_{i-1}^- 
\times M_i$ is an immersed Lagrangian correspondence by Theorem 
\ref{thm:lagr-immers}.  Moreover, since each $M_i$ contains only classes of 
irreducible representations, these Lagrangian correspondences are embedded (cf.  
the proof of Proposition \ref{prop:lagr-embed}).  The symplectic manifolds 
$M_0$ and $M_4$ can each be identified with the single-pointed space $\mathscr 
R_{0,3} = \left\{ [\mathbf i, \mathbf j, -\mathbf k] \right\}$.  Altogether the 
decomposition \eqref{eq:tangle-decomp} induces a cyclic generalized Lagrangian 
correspondence $(L_{01}, L_{12}, L_{23}, L_{34})$ from $M_0 = \left\{ 
\text{pt.} \right\}$ to itself.

By comparing the descriptions of $Y_{01}^*$ and $Y_{34}^*$ above with the proof 
of Lemma \ref{lem:hbdy-xvar}, it is evident that the induced Lagrangian 
correspondences may be expressed as
\begin{equation} \label{eq:l01}
L_{01} = \left\{ (\text{pt.}, [\rho]) \mid [\rho] \in \mathscr R_{g,2k+3}, \; 
\rho([\alpha_i]) = 1 \right\},
\end{equation}
and
\begin{equation} \label{eq:l34}
L_{34} = \left\{ ([\rho], \text{pt.}) \mid [\rho] \in \mathscr R_{g,2k+3}, \; 
\rho([\beta_i]) = 1 \right\},
\end{equation}
where the $\alpha$- and $\beta$-curves are induced by the Morse function in the 
usual manner, and that these may be identified with $L_\alpha$ and $L_\beta$, 
resp.  As we shall see below, the Lagrangians $L_\alpha', L_\beta'$ induced 
after stabilization may be identified with two of the geometric compositions in 
our generalized correspondence, and the embedded composition theorem in quilted 
Floer homology will permit to relate the Floer homologies of these Lagrangian 
pairs.  We will give descriptions of the remaining Lagrangian correspondences 
along with their geometric compositions below in the course of establishing 
these facts.

Recall that the \emph{geometric composition} $L_{(i-1)i} \circ L_{i(i+1)} 
\subset M_{i-1}^- \times M_{i+1}$ is the image under the projection of the 
intersection
\begin{equation} \label{eq:gcomp-inter}
L_{(i-1)i} \times L_{i(i+1)} \cap M_{i-1}^- \times \Delta_{M_i} \times 
M_{i+1},
\end{equation}
where $\Delta_{M_i} \subset M_i^- \times M_i$ is setwise the diagonal.  When 
the intersection \eqref{eq:gcomp-inter} is transversal and its projection to 
$M_{i-1}^- \times M_{i+1}$ injective, the geometric composition is said to be 
\emph{embedded}.  The embedded composition theorem \cite[Theorem 
5.4.1]{WW-qltd_floer_cohom} states that, when the generalized Lagrangian 
  correspondence $(L_{01}, \dotsc, L_{(i-1)i}, L_{i(i+1)}, \dotsc, L_{r(r+1)})$ 
satisfies certain technical conditions and if $L_{(i-1)i} \circ L_{i(i+1)}$ is 
an embedded geometric composition, then the quilted Floer homologies 
$\mathrm{HF}(L_{01}, \dotsc, L_{(i-1)i}, L_{i(i+1)}, \dotsc, L_{r(r+1)})$ and 
$\mathrm{HF}(L_{01}, \dotsc, L_{(i-1)i} \circ L_{i(i+1)}, \dotsc, L_{r(r+1)})$ 
are related by a canonical, grading-preserving isomorphism.  We will apply this 
theorem to exhibit such an isomorphism between $\mathrm{SIK}(\mathcal H) = 
\mathrm{HF}(L_\alpha, L_\beta)$ and $\mathrm{SIK}(\mathcal H') = 
\mathrm{HF}(L_\alpha', L_\beta')$ by characterizing the three geometric 
compisitions arising here according to the following:

\begin{lem}[] \label{lem:gcomp-descr}
The ``outer'' geometric compositions $L_{01} \circ L_{12}$ and $L_{23} \circ 
L_{34}$ can be identified with the Lagrangians $L_\alpha'$ and $L_\beta'$, 
resp.,  in $\mathscr R(\Sigma_2^*)$.  By considering $M_1 \times M_3$ as 
$\mathscr R_{g,2k+3} \times \mathscr R_{g,2k+3}$ the ``inner'' composition 
$L_{12} \circ L_{23}$ is equal to the diagonal $\Delta_{13} \subset M_1^- 
\times M_3$ in the case of either surface or link stabilization.  Moreover, all 
three of these geometric compositions are embedded.
\end{lem}

After showing that $(L_{01}, L_{12}, L_{23}, L_{34})$ meets the aforementioned 
technical conditions, we use Lemma \ref{lem:gcomp-descr} to arrive at the 
following computation on quilted Floer homologies:
\begin{align} \label{eq:emb-comp}
\mathrm{HF}(L_{01}, \Delta_{13}, L_{34})
&= \mathrm{HF}(L_{01}, L_{12} \circ L_{23}, L_{34}) \\
&\cong \mathrm{HF}(L_{01}, L_{12}, L_{23}, L_{34}) \\
\label{eq:2comp} &\cong \mathrm{HF}(L_{01} \circ L_{12}, L_{23} \circ L_{34}).
\end{align}
Since the first of these is equivalent to $\mathrm{HF}(L_\alpha, L_\beta)$ and 
the last to $\mathrm{HF}(L_\alpha', L_\beta')$, we conclude that 
$\mathrm{SIK}(\mathcal H)$ and $\mathrm{SIK}(\mathcal H')$ are canonically 
isomorphic as relatively graded abelian groups, proving Theorem 
\ref{thm:heeg-mv-invar}.

\begin{proof}[Proof of Lemma \ref{lem:gcomp-descr}]
In the following let $R(X^*)$ denote the $\mathrm{SU}(2)$-representation 
variety of the tangle or punctured surface $X^*$ subject to the condition that 
meridians of the implicitly specified codimension-$2$ submanifold be taken to 
traceless elements---prior to taking the quotient by conjugation.  In the 
absence of such a codimension-$2$ submanifold, this is just the usual 
representation variety.

\textbf{Case 1 (surface stabilization):}
In topological terms the punctured surface $\Sigma_2^*$ is the connected sum of 
$\Sigma_1^*$ (or $\Sigma_3^*$) with a torus containing the new attaching curves 
$\alpha', \beta'$.  Let $\gamma$ be the boundary of the disk along which this 
sum is formed, so that $\Sigma_2^*$ decomposes along $\gamma$ into two pieces: 
one the torus minus a disk, which we denote $S$, and the other $\Sigma_2^* 
\setminus S$ diffeomorphic to $\Sigma_1^* \setminus D$, the complement of a 
disk in $\Sigma_1^*$ disjoint from all punctures and attaching curves.  By the 
Seifert--Van Kampen theorem, the traceless representation variety of 
$\Sigma_2^*$ takes the form of the twisted product $R(\Sigma_2^* \setminus S) 
\times_{R(A)} R(S)$, where $A$ is an annular neighborhood of $\gamma$.  Since 
$\pi_1(S)$ is the free group $\langle [\alpha'] \rangle * \langle  [\beta'] 
\rangle$, the representation variety $R(S)$ is just $\mathrm{SU}(2) \times 
\mathrm{SU}(2)$.

The curve $\gamma$ bounds a disk in $Y_{12}^*$ along which the tangle likewise 
splits into two pieces: one diffeomorphic to the cylindrical tangle $\Sigma_1^* 
\times [0,1]$, and the other a solid torus $V$ with $S \subset \partial V$.  We 
therefore have $R(Y_{12}^*) = R(\Sigma_1^* \times [0,1]) \times R(V)$.  Here we 
have $\pi_1(V) = \left\langle [\beta'] \right\rangle$, hence $R(V) = 
\mathrm{SU}(2)$.  The pullback map $R(Y_{12}^*) \to R(\Sigma_2^*)$ respects the 
product structures above, decomposing as the product of pullback maps 
$R(\Sigma_1^* \times [0,1]) \to R(\Sigma_2^* \setminus S)$ and $R(V) \to R(S)$.  
The latter of these is the inclusion of the second factor $\mathrm{SU}(2) 
\hookrightarrow \mathrm{SU}(2) \times \mathrm{SU}(2)$, sending $\rho\vert_V$ to 
its restriction $\rho\vert_S$ satisfying $\rho\vert_S([\alpha']) = 1$ and 
$\rho\vert_S([\beta']) = \rho\vert_V([\beta'])$.  Since the inclusion of 
$\Sigma_2^* \setminus S$ into the cylindrical part $\Sigma_1^* \times [0,1]$ of 
$Y_{12}^*$ induces a surjective map on fundamental groups, the pullback map 
$R(\Sigma_1^* \times [0,1]) \to R(\Sigma_2^* \setminus S)$ is likewise 
injective.  Since any $\rho\vert_{\Sigma_2^* \setminus S} \times \rho\vert_S 
\in \operatorname{im}(R(Y_{12}^*) \to R(\Sigma_2^*))$ is necessarily trivial on 
$A$, it follows that $\rho\vert_{\Sigma_2^* \setminus S}$ here extends uniquely 
over $\Sigma_1^* \times [0,1] \subset Y_{12}^*$.  We therefore conclude that 
$R(Y_{12}^*) \to R(\Sigma_2^*)$ is a smooth embedding, and by 
$\mathrm{SU}(2)$-equivariance, so too is the projection $L_{12} \to M_2$ 
resulting from taking the quotient by conjugation.

Because $\Sigma_1^*$ is a deformation retract of the cylindrical part 
$\Sigma_1^* \times [0,1]$ of $Y_{12}^*$, we may view $R(\Sigma_1^*)$ and 
$R(\Sigma_1^* \times [0,1])$ as being canonically identified up to conjugacy.  
The pullback map $R(Y_{12}^*) \to R(\Sigma_1^*)$ may then be interpreted as the 
projection from $R(Y_{12}^*) = R(\Sigma_1^* \times [0,1]) \times R(V)$ onto the 
first factor, so that $L_{12} \to M_1^-$ is a smooth submersion.  We then see 
that the restrictions $\rho\vert_{\Sigma_2^* \setminus S}$ of traceless 
representations $\rho \in \operatorname{im}(R(Y_{12}^*) \to R(\Sigma_2^*))$ are 
naturally identified with those on $\Sigma_1^*$, and $\rho$ is then uniquely 
determined by the choice of $\rho([\beta']) \in \mathrm{SU}(2)$.  From this it 
is clear that the geometric composition $L_{01} \circ L_{12}$ may be identified 
with $L_{\alpha}'$.  The correspondence $L_{23}$ admits a similar description, 
with $\rho \in \operatorname{im}(R(Y_{23}^*) \to R(\Sigma_2^*))$ characterized 
by the condition $\rho\vert_S([\beta']) = 1$, while $\rho\vert_S([\alpha'])$ 
may be chosen arbitrarily in $\mathrm{SU}(2)$, and it follows from the same 
line of reasoning that $L_{23} \circ L_{34}$ may be identified with $L_\beta'$.

The injectivity of the projection $L_{01} \times_{M_1} L_{12} \to L_{01} \circ 
L_{12}$ is a consequence of the injectivity of $R(Y_{12}^*) \to R(\Sigma_2^*)$ 
and the $\mathrm{SU}(2)$-equivariance of the pullback maps.  If $(\text{pt.}, 
[\rho_2]) \in L_{01} \circ L_{12}$, then $\rho_2 \in R(\Sigma_2^*)$ extends 
uniquely over $Y_{12}^*$, and this extension in turn corresponds uniquely with 
some $\rho_1 \in R(\Sigma_1^*)$ determined by $\rho_2\vert_{\Sigma_2^* 
\setminus S}$.  Hence $(\text{pt.}, [\rho_1], [\rho_1], [\rho_2]) \in M_0 
\times \Delta_{M_1} \times M_2$ is the unique element in the preimage of the 
projection of $(\text{pt.}, [\rho_2])$.

That $L_{01} \times_{M_1} L_{12}$ is cut out transversally by the diagonal 
$\Delta_{M_1}$ follows from $\mathrm{SU}(2)$-equivariance and the fact that 
$R(Y_{12}^*) \to R(\Sigma_1^*)$ is surjective, hence $L_{12} \to M_1^-$ a 
smooth submersion.  Take an arbitrary tangent vector $(0, \xi_1, \xi_1', \xi_2) 
\in T_{[\rho]}(\text{pt.} \times M_1 \times M_1 \times M_2)$.  There can be 
found a vector of the form $(0, 0, \xi_1' - \xi_1, \xi_2') \in 
T_{[\rho]}(L_{01} \times L_{12})$.  Then summing with $(0, \xi_1, \xi_1, \xi_2 
- \xi_2') \in T_{[\rho]}(M_0 \times \Delta_1 \times M_2)$ gives the original 
tangent vector, and we conclude that the intersection is transverse.  Therefore 
$L_{01} \circ L_{12}$ is embedded.  The arguments for the embeddedness of 
$L_{23} \circ L_{34}$ are essentially the same.

Finally, we consider the composition $L_{12} \circ L_{23}$.  Representations 
$\rho_1 \in \operatorname{im}(R(Y_{12}^*) \to R(\Sigma_2^*))$ and $\rho_2 \in 
\operatorname{im}(R(Y_{23}^*) \to R(\Sigma_2^*))$ agree up to conjugation if 
and only if their restrictions to $\Sigma_2^* \setminus S$ are conjugate and 
are both trivial along $\gamma$.  Two such representations correspond uniquely 
with represenations in $R(\Sigma_1^*)$ and $R(\Sigma_3^*)$, resp., and these 
necessarily represent the same conjugacy class in $\mathscr R_{g,2k+3}$.  From 
this we conclude that the geometric composition $L_{12} \circ L_{23}$ is 
set-theoretically the diagonal $\Delta_{13} \subset M_1 \times M_3 = \mathscr 
R_{g,2k+3} \times \mathscr R_{g,2k+3}$.  The uniqueness further implies that 
the projection from $L_{12} \times_{M_2} L_{23}$ is injective.

The transversality of $L_{12} \times L_{23} \cap M_1 \times \Delta_{M_2} \times 
M_3$ can be deduced from the fact that the images of $R(Y_{12}^*)$ and 
$R(Y_{23}^*)$ intersect transversally in $R(\Sigma_2^*)$, as can be worked out 
simply in terms of the dimensions of the representation varieties.  This 
clearly implies that the intersection of $\operatorname{im}(R(Y_{12}^*)) \times 
\operatorname{im}(R(Y_{23}^*)) \subset R(\Sigma_2^*) \times R(\Sigma_2^*)$ with 
the diagonal is transversal, and then transversality at the level of character 
varieties follows, since the quotient map induced by the conjugation action is 
a smooth submersion.

\textbf{Case 2 (link stabilization):}
Recall that we now have $\alpha', \beta'$ bounding disks $A', B' \subset 
\Sigma_2^*$, resp., with two new punctures $w' \in A' \cap B'$ and $z' \in B'$ 
and exactly one of the original punctures $z_0 \in B'$ as well.  Take a smooth 
curve $\gamma \subset \Sigma_2^*$ which is a slight pushoff of $\partial(A' 
\cup B')$.  The gradient flow of $\gamma$ is diffeomorphic to the cylinder $S^1 
\times [1,2]$ which splits $Y_{12}^*$ into two pieces---one cylindrical and one 
non-cylindrical tangle.  The non-cylindrical component is equivalent to a disk 
times an interval, containing two arcs, which we write altogether as $(D \times 
[1, 2])^*$.  One end $D_1^* = (D \times \left\{ 1 \right\})^* \subset 
\Sigma_1^*$ is a disk with a single puncture identified with $z_0$.  The other 
$D_2^* = (D \times \left\{ 2 \right\})^* \subset \Sigma_2^*$ is a disk with the 
three punctures $z_0, w', z'$.  Within $(D \times [1, 2])^*$ is one arc 
connecting $z_0 \times \left\{ 1 \right\}$ with $z_0 \times \left\{ 2 \right\}$ 
and a cup connecting $w'$ with $z'$.  The gradient flow of the Morse function 
serves to identify the cylindrical component with $\Sigma_g^* \setminus D 
\times [1, 2]$, where $D$ is a disk and the punctures in $\Sigma_g^*$ are 
identified with $\mathbf w \cup \mathbf z \setminus \left\{ z_0 \right\} \cup 
\mathbf q$, so that $\Sigma_g^* \setminus D \times \left\{ i \right\} = 
\Sigma_i^* \setminus D_i^*$ for $i = 1, 2$.  By the Seifert--Van Kampen theorem 
we obtain the following expressions of the representation varieties as twisted 
products:
\begin{equation} \label{eq:surf-repvar2}
R(\Sigma_2^*) = R(\Sigma_g^* \setminus D) \times_{R(S^1)} R(D_2^*),
\end{equation}
\begin{equation} \label{eq:tangle-repvar12}
R(Y_{12}^*) = R(\Sigma_g^* \setminus D \times [1, 2]) \times_{R(S^1 \times [1, 
2])} R((D \times [1, 2])^*),
\end{equation}
\begin{equation} \label{eq:surf-repvar1}
R(\Sigma_1^*) = R(\Sigma_g^* \setminus D) \times_{R(S^1)} R(D_1^*).
\end{equation}
Note with regard to $R(\Sigma_1^*)$ that the restriction to $\partial D$ 
completely determines the restriction to $D_1^*$, hence $R(\Sigma_1^*)$ may be 
identified with its first factor $R(\Sigma_g^* \setminus D)$.

The pullback maps $R(Y_{12}^*) \to R(\Sigma_1^*), R(\Sigma_2^*)$ again respect 
the product structures.  The restrictions to the first factor are induced by 
the inclusion of a deformation retract, hence we may again think of each of the 
first factors in the products above as being related by the identity map, and 
$R(Y_{12}^*) \to R(\Sigma_1^*)$ as the projection on the first factor.  Since 
the map $\pi_1(D_2^*) \to \pi_1((D \times [1, 2])^*)$ is surjective, the 
pullback map $R(Y_{12}^*) \to R(\Sigma_2^*)$ is a smooth embedding.  Let $c_0, 
c_1, c_2 \in \pi_1(D_2^*)$ be the classes of standard meridians about $z_0, w', 
z'$, resp.  By writing representations in $\operatorname{im}(R(Y_{12}^*) \to 
R(\Sigma_2^*))$ in the form $\rho\vert_{\Sigma_g^* \setminus D} \times_{\gamma} 
\rho\vert_{D_2^*}$, we see that the image is cut out by the equation 
$\rho\vert_{D_2^*}(c_1) = -\rho\vert_{D_2^*}(c_2)$.  This is equivalent to the 
condition that $\rho\vert_{D_2^*}([\alpha']) = 1$.  The image of $c_0$ on the 
other hand is arbitrary, and this together with $\rho\vert_{\Sigma_g^* 
\setminus D}$ uniquely determine a representation in $R(\Sigma_1^*)$.  Hence we 
conclude as above that $L_{01} \circ L_{12}$ may be identified with 
$L_{\alpha}'$.  The embeddedness of this composition follows from precisely the 
same argument as in the case of surface stabilization above, and again 
analogous descriptions of the representation varieties lead to the same 
conclusions for $L_{23} \circ L_{34}$.

The image of $R(Y_{23}^*) \to R(\Sigma_2^*)$ is cut out by the equation 
$\rho\vert_{D_2^*}(c_0) = -\rho\vert_{D_2^*}(c_1)$.  Representations $\rho_1 
\in \operatorname{im}(R(Y_{12}^*) \to R(\Sigma_2^*))$ and $\rho_2 \in 
\operatorname{im}(R(Y_{23}^*) \to R(\Sigma_2^*))$ are conjugate if and only if 
both satisfy $\rho_i(c_0) = -\rho_i(c_1) = \rho_i(c_2)$.  Similarly as in the 
previous case, the account above shows that these representations correspond 
uniquely with ones on $R(\Sigma_1^*)$ and $R(\Sigma_3^*)$, resp., which 
necessarily represent the same class in $\mathscr R_{g,2k+3}$.  Hence, after 
passing to the quotient, we again have that $L_{12} \circ L_{23}$ is equal to 
the diagonal $\Delta_{13}$ on $\mathscr R_{g,2k+3}$, with the projection from 
$L_{12} \times_{M_2} L_{23}$ being injective.  Transversality follows by the 
same argument as given for surface stabilization.
\end{proof}

To establish the validity of our application of the embedded composition 
theorem above, we now check the outstanding technical conditions regarding the 
Lagrangian correspondences' monotonicity, relative spin structures and disk 
numbers.

\textbf{Monotonicity:}
The embedded composition theorem requires that all of the symplectic manifolds 
$M_i$, for $i = 0, \dotsc, 3$, be compact and monotone with the same 
monotonicity constant.  The compactness condition is obviously satisfied, and 
that of monotonicity is established in Proposition \ref{prop:surf-monot}.  The 
theorem further requires that each of the Lagrangian correspondences 
$L_{i(i+1)}$ be a monotone Lagrangian submanifold of minimal Maslov number $N_L 
\ge 2$.  We have already seen this verified for $L_{01} = L_\alpha$ and $L_{34} 
= L_\beta$.  The same holds for $L_{12}$ and $L_{23}$ by Proposition 
\ref{prop:sc-monot} and the subsequent remarks since both of these 
correspondences are likewise simply connected.  From the descriptions in the 
proof of Lemma \ref{lem:gcomp-descr} it can be seen that, for instance, 
$R(Y_{12}^*)$ is diffeomorphic to $R(\Sigma_1^*) \times \mathrm{SU}(2)$ in the 
case of surface stabilization and to $R(\Sigma_1^*) \times \mathcal C_0$ in 
that of link stabilization.  Since all representations in $R(\Sigma_1^*)$ have 
central stabilizer, the quotient by conjugation $L_{12}$ has the structure of a 
fiber bundle over the simply connected space $M_1 = \mathscr R_{g,2k+3}$ with 
fiber $\mathrm{SU}(2)$ or $\mathcal C_0$, and is therefore simply connected.  
An analogous argument holds for $L_{23}$.

\textbf{Relative spin structures:}
What is required here is that, for $i = 0, \dotsc, 3$, there should be classes 
$b_i \in H^2(M_i; \mathbb Z/2)$ and $b_{i+1} \in H^2(M_{i+1}; \mathbb Z/2)$ 
such that the sum of pullbacks under the projections $\pi_i^*b_i + 
\pi_{i+1}^*b_{i+1} \in H^2(M_i \times M_{i+1}; \mathbb Z/2)$ constitutes a 
background class for a relative spin structure on $L_{i(i+1)}$.  In \cite[Lemma 
3.19]{WW-floer_fld_thry_tngl} Wehrheim and Woodward exhibit such background 
  classes for Lagrangian correspondences associated to elementary tangles, in 
the general case of a compact Lie group $G$ and suitable choices of conjugacy 
classes.  Their argument can be generalized to non-elementary tangles, but due 
to our choice of $G = \mathrm{SU}(2)$ and the conjugacy class $\mathcal C_0$ of 
traceless elements, we require only what amounts to a special case of Wehrheim 
and Woodward's construction.  Namely, we claim that all of our Lagrangian 
correspondences are embeddings of spin manifolds---as we have already seen of 
$L_{01} = L_\alpha$ and $L_{34} = L_\beta$---so that it suffices to take $b_i = 
0$ for all $i$.

To see that $L_{12}$ and $L_{23}$ are spin, we start by obvserving that, for 
instance, $R(Y_{12}^*)$ is a regular level set of a smooth submersion $\Phi : X 
\to \mathrm{SU}(2)$ where $X$ is a cartesian product of copies of 
$\mathrm{SU}(2)$ and $\mathcal C_0$.  The normal bundle of a regular level set 
is trivial, being isomorphic to the pullback of the tangent space at a point.  
It follows that the second Stiefel-Whitney class $w_2(R(Y_{12}^*))$ is given by 
the restriction $w_2(X)\vert_{R(Y_{12}^*)}$.  Since $\mathrm{SU}(2)$ is 
equivariantly spin, the product $X$ has a unique $\mathrm{SU}(2)$-equivariant 
spin structure, and the same must then be true of $R(Y_{12}^*)$.  After passing 
to the quotient by conjugation, we obtain a spin structure on $L_{12}$, and 
similarly for $L_{23}$.

\textbf{Disk numbers:}  The identifications of $L_{01} \circ L_{12}$ and 
$L_{23} \circ L_{34}$ with $L_{\alpha}'$ and $L_{\beta}'$, resp., imply that 
$\Phi(L_{01} \circ L_{12}) = \Phi(L_{\alpha}') = \Phi(L_{\beta}') = - 
\Phi(L_{23} \circ L_{34})$.  The sign change in the right-most term is a result 
of the composition $L_{23} \circ L_{34}$ residing in $M_2^- \times M_0$ with 
the sign of the symplectic form switched on $M_2$.  Consequently the condition 
on the quilted Floer homology in \eqref{eq:2comp} that $\Phi(L_{01} \circ 
L_{12}) + \Phi(L_{23} \circ L_{34}) = 0$ is satisfied.  As noted in the proof 
of \cite[Theorem 
5.4.1]{WW-qltd_floer_cohom}, it is enough to verify that the disk numbers of 
the Lagrangian correspondences sum to $0$ for just one of the generalized 
correspondences appearing in \eqref{eq:emb-comp}.  That this holds for the 
others may then be deduced \textit{a posteriori} from the embedded composition 
theorem.

This completes the proof of Theorem \ref{thm:heeg-mv-invar}.

\subsection{Naturality}

It would be preferable to have $\mathrm{SIK}(Y, K)$ assign a concrete group to 
the pair $(Y, K)$, rather than an isomorphism class of groups, i.e. to be a 
\emph{natural} homology invariant.  Only for natural invariants does it make 
sense to speak of e.g. isomomorphisms of the Floer homology induced by 
cobordisms.

Loosely speaking, given a pair of Heegaard diagrams $\mathcal H$ and $\mathcal 
H'$ as above, two different sequences of Heegaard moves $\mathcal H 
\longrightarrow \mathcal H'$ correspond to a loop $\mathcal H_t$ in the space 
of compatible diagrams based at $\mathcal H$, and the obstruction to naturality 
is the monodromy of $\mathrm{SIK}$ along $\mathcal H_t$.  Juh\'asz, Thurston 
and Zemke \cite{JTZ-natural} furnish a set of four conditions which being met 
imply that an algebraic object associated to pointed Heegaard diagrams has 
trivial monodromy and is therefore a natural invariant.  Horton uses this with 
an application of quilted Floer homology, along the lines of the one used 
above, to demonstrate naturality for the $3$-manifold symplectic instanton 
homology $\mathrm{SI}(Y)$, and so it would seem that the same should hold here:

\begin{conj}[Naturality]
As a Lagrangian Floer homology invariant of a multi-pointed Heegaard diagram,  
$\mathrm{SIK}(\mathcal H)$ satisfies the four conditions of a \emph{strong 
Heegaard invariant}, see \cite[Definition 2.32]{JTZ-natural}.  There is 
therefore a canonical group representing its isomorphism class, which one may 
take to be the definition of $\mathrm{SIK}(Y, K)$.
\end{conj}

%% file: sec5-examples.tex
\section{Examples and computations} \label{sec:examples}

\begin{exm}[Empty link] \label{exm:empty}
If $K = \emptyset$ is the empty link, then a Heegaard diagram compatible with 
$(Y, \emptyset)$ is just a triply-pointed diagram of the kind used to construct 
Horton's symplectic instanton homology $\mathrm{SI(Y)}$ of the $3$-manifold.  
We may therefore simply define $\mathrm{SIK}(Y, \emptyset) := \mathrm{SI}(Y)$, 
which is a concrete homology group owing to the proven naturality of the latter 
theory \cite[Cor. 5.6]{Hor-sih_trless_xvar}.  In particular $\mathrm{SIK}(S^3, 
\emptyset) \cong \mathbb Z$ and $\mathrm{SIK}(S^2 \times S^1, \emptyset) \cong 
H^*(S^3) \cong \mathbb Z \oplus \mathbb Z$ (see \cite[Sec.  
9]{Hor-sih_trless_xvar}).
\end{exm}

\begin{exm}[Unlinks in $S^3$] \label{exm:unlink}
Let $Y = S^3$ and $K = U^\ell$ be the unlink with $\ell \ge 1$ components.  For 
this we may choose a genus-$0$ Heegaard diagram $\mathcal H_0$ with respect to 
which each component of $U^\ell$ is in $1$-bridge position.  In this case it is 
easy to see that the induced Lagrangians $L_\alpha = L_\beta = L \approx 
(S^2)^{\times \ell}$ coincide in $\mathscr R_{0,2\ell+3}$, hence 
$\mathrm{SIK}(\mathcal H_0) = \mathrm{HF}(L)$ is the self-Floer homology of 
$L$.

Let $f : L \to \mathbb R$ be a $C^2$-small Morse function.  This induces a 
Hamiltonian perturbation $L'$ of $L$ lying in a Weinstein neighborhood of the 
latter, with $\text{Crit}(f) = L \cap L'$.  We may arrange that $f$ has a 
unique local minimum $x_0$ and take this as the base intersection point in 
grading $0 \!\!\mod \mathbb Z/2$, thereby fixing the relative $\mathbb 
Z/2$-grading on the self-Floer cochain complex $\mathrm{CF}^*(L) = 
\mathrm{CF}^*(L, L')$.  In exhibiting his spectral sequence relating self-Floer 
cohomology to singular cohomology, Oh \cite{Oh-spec} demonstrates that 
$\mathrm{CF}^*(L)$ can be identfied with the Morse cochain complex $C_f^*(L)$ 
in such a way that the $\mathbb Z/N_L$-grading on the Floer complex aligns with 
usual grading of the Morse complex:
\begin{equation} \label{eq:floer-morse}
\mathrm{CF}^{\, i \!\!\!\!\mod N_L}(L) = \bigoplus_{\, j \equiv i \!\!\!\!\mod 
N_L} C_f^j(L).
\end{equation}
Since both the Floer and Morse coboundary maps raise degree by $1$, and since 
$N_L$ is necessarily even, it is evident from \eqref{eq:floer-morse} that, if 
$L$ admits a Morse function with all critical points residing in even degree, 
the self-Floer cohomology is isomorphic to the Morse cohomology, and thus to 
singular cohomology.  Hence we obtain the self-Floer homology, which is the 
dual of and isomorphic to the Floer cohomology.

Because in our case $L$ is a product of $2$-spheres, we may choose $f$ to be a 
perfect Morse function, with all critical points of even index.  Therefore 
$\mathrm{SIK}(\mathcal H_0) = \mathrm{HF}(L) \cong H^*((S^2)^{\times \ell})$.
\end{exm}

%% file: sik.bbl
\providecommand{\bysame}{\leavevmode\hbox to3em{\hrulefill}\thinspace}
\providecommand{\MR}{\relax\ifhmode\unskip\space\fi MR }
\providecommand{\MRhref}[2]{%
  \href{http://www.ams.org/mathscinet-getitem?mr=#1}{#2}
}
\providecommand{\href}[2]{#2}